\documentclass{amsart}

\usepackage{paralist}
\usepackage{url}

\usepackage{amssymb,amsmath,amsfonts}
\usepackage{algorithm,algpseudocode}

\usepackage{multirow}
\usepackage[table]{xcolor}
\usepackage{slashbox}
\usepackage{diagbox}
\usepackage{booktabs}

\usepackage{graphicx,here}
\usepackage{subfig}
\usepackage{epsfig}

\theoremstyle{definition}

\title[MO approach to optimal control for a dengue transmission model]{Multiobjective
approach to optimal control\\ for a dengue transmission model}

\thanks{This is a preprint
of a paper whose final and definite form is in
\emph{Statistics, Optimization \& Information Computing} (SOIC),
ISSN 2310-5070 (online), ISSN 2311-004X (print). Paper submitted 28/April/2015;
accepted, after a revision, 25/June/2015.}

\author[Denysiuk, Rodrigues, Monteiro, Costa, Esp\'{i}rito Santo and Torres]{%
Roman Denysiuk$^{1}$\\
{\tt \lowercase{roman.denysiuk@algoritmi.uminho.pt}}\vspace*{0.2cm}\\
Helena Sofia Rodrigues$^{2, 3}$\\
{\tt \lowercase{sofiarodrigues@esce.ipvc.pt}}\vspace*{0.2cm}\\
M. Teresa T. Monteiro$^{1, 4}$\\
{\tt \lowercase{tm@dps.uminho.pt}}\vspace*{0.2cm}\\
Lino Costa$^{1, 4}$\\
{\tt \lowercase{lac@dps.uminho.pt}}\vspace*{0.2cm}\\
Isabel Esp\'{i}rito Santo$^{1, 4}$\\
{\tt \lowercase{iapinho@dps.uminho.pt}}\vspace*{0.2cm}\\
Delfim F. M. Torres$^{3, 5}$\\
{\tt \lowercase{delfim@ua.pt}}}



\begin{document}

\maketitle

\centerline{\footnotesize
$^1$Algoritmi R\&D Center, University of Minho, 4710-057 Braga, Portugal}

\medskip

\centerline{\footnotesize
$^2$School of Business Studies, Viana do Castelo Polytechnic Institute,
4930-678 Valen\c{c}a, Portugal}

\medskip

\centerline{\footnotesize
$^3$CIDMA R\&D Center, University of Aveiro, 3810-193 Aveiro, Portugal}

\medskip

\centerline{\footnotesize
\footnotesize $^4$Department of Production and Systems Engineering,
University of Minho, 4710-057 Braga, Portugal}

\medskip

\centerline{\footnotesize
$^5$Department of Mathematics, University of Aveiro, 3810-193 Aveiro, Portugal}

\bigskip


\begin{abstract}
During the last decades, the global prevalence of dengue progressed dramatically.
It is a disease which is now endemic in more than one hundred countries of Africa,
America, Asia and the Western Pacific. This study addresses a mathematical model
for the dengue disease transmission and finding the most effective ways
of controlling the disease. The model is described by a system of ordinary
differential equations representing human and vector dynamics. Multiobjective
optimization is applied to find the optimal control strategies, considering
the simultaneous minimization of infected humans and costs due to insecticide
application. The obtained results show that multiobjective optimization is an
effective tool for finding the optimal control. The set of trade-off solutions
encompasses a whole range of optimal scenarios, providing valuable information
about the dynamics of infection transmissions. The results are discussed
for different values of model parameters.
\end{abstract}


\section{Introduction}
\label{intr}

Dengue is a vector-borne disease transmitted from an infected human to a female
\emph{Aedes} mosquito by a bite. The mosquito, which needs regular meals of blood
to feed their eggs, bites a potentially healthy human and transmits the disease,
turning it into a cycle. There are four distinct, but closely related, viruses
that cause dengue. The four serotypes, named DEN-1 to DEN-4, belong to the
\emph{Flavivirus} family, but they are antigenically distinct. Recovery from
infection by one serotype provides lifelong immunity against that serotype but
provides only partial and transient protection against subsequent infection
by the other three viruses. There are strong evidences that a sequential
infection increases the risk of developing dengue hemorrhagic fever.

The spread of dengue is attributed to the geographic expansion of the mosquitoes
responsible for the disease: \emph{Aedes aegypti} and \emph{Aedes albopictus}.
The \emph{Aedes aegypti} mosquito is a tropical and subtropical species widely
distributed around the world, mostly between latitudes $35^{o}$N and 35$^o$S.
In urban areas, \emph{Aedes} mosquitoes breed on water collections in artificial
containers such as cans, plastic cups, used tires, broken bottles and flower pots.
Due to its high interaction with humans and its urban behavior,
the \emph{Aedes aegypti} mosquito is considered the major responsible
for the dengue transmission. The life cycle of a mosquito has four distinct stages:
egg, larva, pupa and adult. In the case of \emph{Aedes aegypti}, the first three
stages take place in, or near, the water, whereas the air is the medium
for the adult stage~\cite{Otero2008}.

It is very difficult to control or eliminate \emph{Aedes aegypti} mosquitoes
due to their resiliency, fast adaptation to changes in the environment
and their ability to rapidly bounce back to initial numbers after disturbances
resulting from natural phenomena (e.g., droughts) or human interventions
(e.g., control measures). Primary prevention of dengue resides mainly
in mosquito control. There are two primary methods: larval control
and adult mosquito control, depending on the intended target.
\emph{Larvicide} treatment is done through long-lasting chemical in order
to kill larvae and preferably have World Health Organization clearance
for use in drinking water~\cite{Derouich2003}. \emph{Adulticides}
is the most common measure. Its application can have a powerful impact
on the abundance of adult mosquito vector. However, the efficacy is often
constrained by the difficulty in achieving sufficiently high coverage
of resting surfaces~\cite{Devine2009}.

The present study addresses a problem of finding the most effective ways
of controlling the dengue disease. To this end, a mathematical model
for the dengue transmission, including a control variable presented
in \cite{RodriguesMonteiro2012}, is adopted. The main contributions are:
(i) adapting multiobjective optimization methodologies for finding the
optimal control in a mathematical model for the dengue transmission,
(ii) analysis of the search ability of different methods on the resulting
optimization problem, (iii) discussion of results for the optimal control
problem, emphasizing advantages of the proposed approach when compared
with results available in the literature. The aim of the work is to promote
multiobjective optimization in epidemiological research community and practice.

The remainder of this paper is organized as follows. Section~\ref{sec:model}
describes a mathematical model for the dengue disease transmission. The problem
of finding the optimal control is formulated in Section~\ref{sec:problem},
including traditional and proposed approaches. Section~\ref{sec:mo} presents
some general concepts in multiobjective optimization and methods adopted
for solving the problem. Section~\ref{sec:experiments} presents and discusses
our results, obtained by numerical simulations, which includes the comparison
of different multiobjective optimization methods, control strategies and
dengue dynamics that correspond to different optimal scenarios. The study
is ended with Section~\ref{sec:cons} of conclusions
and some directions of future work.


\section{The mathematical model}
\label{sec:model}

This section introduces a mathematical model for the dengue disease transmission.
The model is based on a system of ordinary differential equations, and includes
the real data of a dengue disease outbreak that occurred in
the Cape Verde archipelago in 2009~\cite{ministerio,RMT2013a}.

The model consists of eight mutually-exclusive compartments representing
the human and vector dynamics. It also includes a control parameter,
an adulticide spray, as a measure to fight the disease.

The notation used in the mathematical model includes
four epidemiological states for humans:
\begin{quote}
\begin{tabular}{lcl}
$S_h(t)$ & -- & number of susceptible individuals at time $t$;\\
$E_h(t)$ & -- & number of exposed individuals at time $t$;\\
$I_h(t)$ & -- & number of infected individuals at time $t$; \\
$R_h(t)$ & -- & number of resistant individuals at time $t$.
\end{tabular}
\end{quote}
It is assumed that the total human population
$(N_h)$ is constant, so
$$
N_h=S_h(t)+E_h(t)+I_h(t)+R_h(t).
$$
There are also four other state variables related to the mosquitoes (disease vectors):
\begin{quote}
\begin{tabular}{lcl}
$A_m(t)$ & -- & number of vectors in aquatic phase at time $t$;\\
$S_m(t)$ & -- & number of susceptible vectors at time $t$;\\
$E_m(t)$ & -- & number of exposed vectors at time $t$;\\
$I_m(t)$ & -- & number of infected vectors at time $t$.\\
\end{tabular}
\end{quote}
Similarly, it is assumed that the total adult
mosquito population ($N_m$) is constant,
which means that $N_m=S_m(t)+E_m(t)+I_m(t)$.

The model includes a control variable, which represents the amount of insecticide
that is continuously applied during a considered period, as a measure to fight the disease:
\begin{quote}
\begin{tabular}{lcl}
$c(t)$ & -- & level of insecticide campaigns at time $t$.\\
\end{tabular}
\end{quote}
The control variable $c(t)$ is an adimensional value that
is considered in relative terms, varying from 0 to 1.

In the following, for the sake of simplicity, the independent variable $t$
is omitted when writing the dependent variables
(for instance, $S_h$ is used instead of $S_h (t)$).

The parameters necessary to completely describe the model are presented
in Table~\ref{tab:param}. This set of parameters includes the real data
related to the outbreak of dengue disease occurred
in Cape Verde in 2009 \cite{RMT2013a}.

\begin{table}[h!]
\centering
\begin{tabular}{!{\vrule width 1pt} c|c|c !{\vrule width 1pt}}
\specialrule{1pt}{0pt}{0pt}
Parameter & Description & Value \\ \hline
$N_h$ & total population & 480000 \\[0.2cm]
$B$ & average daily bites (per mosquito per day)& 1 \\[0.2cm]
$\beta_{mh}$ & transmission probability from $I_m$ (per bite)& 0.375  \\[0.2cm]
$\beta_{hm}$ & transmission probability from $I_h$ (per bite)& 0.375  \\[0.2cm]
$1/\mu_{h}$ & average human lifespan (in days)& $1/(71\times365)$  \\[0.2cm]
$1/\eta_{h}$ & mean viremic period (in days)& 1/3 \\[0.2cm]
$1/\mu_{m}$ & average lifespan of adult mosquitoes (in days)& 1/11  \\[0.2cm]
$\varphi$ & number of eggs at each deposit per capita (per day)& 6  \\[0.2cm]
$\mu_{A}$ & natural mortality of larvae (per day)& 1/4  \\[0.2cm]
$\eta_A$ & rate of maturing from larvae to adult (per day)& 0.08  \\[0.2cm]
$1/\eta_m$ & extrinsic incubation period (in days)& 1/11   \\[0.2cm]
$1/\nu_h$ & intrinsic incubation period (in days)& 1/4  \\[0.2cm]
$m$ & number of female mosquitoes per human & 6 \\[0.2cm]
$k$ & number of larvae per human & 3\\
\specialrule{1pt}{0pt}{0pt}
\end{tabular}
\caption{Model parameters.}\label{tab:param}
\end{table}

Furthermore, in order to obtain a numerically stable problem,
all the state variables are normalized as follows:
\begin{center}
\begin{tabular}{llll}
$\displaystyle s_h=\frac{S_h}{N_h},$ & $\displaystyle e_h
=\frac{E_h}{N_h},$ & $\displaystyle i_h=\frac{I_h}{N_h},$
& $\displaystyle r_h=\frac{R_h}{N_h},$\\
& & & \\
$\displaystyle a_m=\frac{A_h}{kN_h},$ & $\displaystyle s_m=\frac{S_m}{mN_h},$
& $\displaystyle e_m=\frac{E_m}{mN_h},$ & $\displaystyle i_m=\frac{I_m}{mN_h}.$\\
\end{tabular}
\end{center}
Thus, the dengue epidemic is modeled by the following nonlinear
time-varying state equations~\cite{RodriguesMonteiro2012}:
\begin{equation}
\label{dengue:model}
\begin{tabular}{l}
$\left\{
\begin{array}{l}
\displaystyle\frac{ds_h}{dt} = \mu_h - (B\beta_{mh}m i_m +\mu_h) s_h\\[0.25cm]
\displaystyle\frac{de_h}{dt} = B\beta_{mh}m i_m s_h - (\nu_h + \mu_h )e_h\\[0.25cm]
\displaystyle\frac{di_h}{dt} = \nu_h e_h -(\eta_h  +\mu_h) i_h\\[0.25cm]
\displaystyle\frac{dr_h}{dt} = \eta_h i_h - \mu_h r_h\\[0.25cm]
\displaystyle\frac{da_m}{dt} = \varphi \frac{m}{k}(1-a_m)(s_m+e_m+i_m)-(\eta_A+\mu_A) a_m\\[0.25cm]
\displaystyle\frac{ds_m}{dt} = \eta_A \frac{k}{m}a_m-(B \beta_{hm}i_h+\mu_m) s_m-c s_m\\[0.25cm]
\displaystyle\frac{de_m}{dt} = B \beta_{hm}i_h s_m-(\mu_m + \eta_m) e_m-c e_m\\[0.25cm]
\displaystyle\frac{di_m}{dt} = \eta_m e_m -\mu_m i_m - c i_m\\[0.25cm]
\end{array}
\right. $\\
\end{tabular}
\end{equation}
subject to the initial conditions
\begin{equation}
\label{dengue:initial}
\begin{array}{llll}
s_h(0)=0.99865, & e_h(0)=0.00035, & i_h(0)=0.001, & r_h(0)=0, \\
a_m(0)=1, & s_{m}(0)=1, & e_m(0)=0, & i_m(0)=0.\\
\end{array}
\end{equation}

Since any mathematical model is an abstraction of a complex natural system,
additional assumptions are made to make the model mathematically treatable.
This is also the case for the above epidemiological model,
which comprises the following assumptions:
\begin{itemize}
\item the total human population ($N_h$) is constant;
\item there is no immigration of infected individuals into the human population;
\item the population is homogeneous, which means that every
individual of a compartment is homogeneously mixed with other individuals;
\item the coefficient of transmission of the disease is fixed
and does not vary seasonally;
\item both human and mosquitoes are assumed to be born susceptible, i.e.,
there is no natural protection;
\item there is no resistant phase for mosquitoes, due to their short lifetime.
\end{itemize}


\section{Problem formulation}
\label{sec:problem}

In this section, the problem of finding the optimal control in the mathematical model
for the dengue transmission is formulated. The optimal control represents the most
effective way of controlling the disease that can be adopted by authorities in
response to its outbreak. First, a traditional approach based on optimal control
theory is presented. Then, the proposed approach, based
on multiobjective optimization, is introduced.


\subsection{Optimal control theory}
\label{sec:oc:app}

A problem of finding a control law for a given system is commonly formulated
and solved using optimal control theory. Our aim is to find the optimal value
$c^*$ of the control $c$, such that the associated state trajectories
$s_h^*$, $e_h^*$, $i_h^*$, $r_h^*$, $a_m^*$, $s_m^*$, $e_m^*$, $i_m^*$
are solution of the system \eqref{dengue:model} in the time interval $[0, T]$,
subject to the initial conditions~\eqref{dengue:initial},
and minimize an objective functional.

Consider the state system of ordinary differential equations~\eqref{dengue:model}
and the set of admissible control functions given by:
\[
\Omega = \{ c(\cdot) \in L^{\infty}(0,T) \, | \,
0 \leq c(t) \leq 1, \forall t \in [0,T] \}.
\]
The objective functional is defined by~\cite{RodriguesMonteiro2013}:
\begin{equation}
\label{functional}
J(c(\cdot)) = \int_{0}^{T}\left[ \gamma_D i_h(t)^2 + \gamma_S c(t)^2\right]dt,
\end{equation}
where $\gamma_D$ and $\gamma_S$ are positive constants representing the costs
weights of infected individuals and spraying campaigns, respectively.
The objective functional given by~\eqref{functional} is a function of state
and control variables, representing cumulative costs due to infected population
and prevention measures. The optimal control problem consists in determining
$\left(s_h^*, \, e_h^*, \, i_h^*, \, r_h^*,\, a_m^*,\, s_m^*,\, e_m^*,
\, i_m^*\right)$, associated to an admissible control $c^*(\cdot) \in \Omega$
on the time interval $[0, T]$, satisfying~\eqref{dengue:model}, the initial
conditions~\eqref{dengue:initial}, and minimizing
the cost functional \eqref{functional}, \textrm{i.e.},
\[
J(c^*(\cdot)) = \min_{c(\cdot) \in \Omega} J(c(\cdot)) \, .
\]
The optimal control can be derived using
Pontryagin's maximum principle~\cite{pontryagin62}.


\subsection{Multiobjective approach}
\label{sec:mo:app}

An approach based on optimal control theory allows to obtain
a single optimal solution. The obtained solution represents
some decision maker's perspective on controlling the disease.
This is reflected by the constants $\gamma_D$ and $\gamma_S$
that must be provided in advance. Though the choice of proper
values for the constants is not an easy task for all the cases.
A single optimal solution provides limited information for
the decision maker, while leaving a large range of alternatives
unexplored. To overcome these drawbacks, the present study decomposes
the cost functional~\eqref{functional} into two separate objectives
and seeks the optimal control, minimizing simultaneously
the costs due to the infected population and the costs
associated with insecticide. The resulting optimization problem
is defined as follows:
\begin{equation}
\label{dengue:mop}
\begin{array}{rl}
\text{minimize:}
& f_1(i_h(\cdot)) = \int_{0}^{T} i_h(t) \, dt \\
& f_2(c(\cdot)) = \int_{0}^{T} c(t) \, dt \\[0.3 cm]
\text{subject to:} & \eqref{dengue:model} \text{ and } \eqref{dengue:initial} \\
\end{array}
\end{equation}
where $T$ is a given period of time, $f_1$ and $f_2$ represent the cost incurred
in the form of infected population and the cost of applying
insecticide for the period $T$, respectively.


\section{Multiobjective optimization}
\label{sec:mo}

This section presents general concepts in multiobjective optimization
and discusses some popular methods for multiobjective optimization based
on scalarization. The outline of an algorithm for finding
optimal solutions is presented. For a method that gives a direct solution
to the full Pareto set, via solutions of a certain PDE, we refer the reader
to \cite{MR1198561,MR1447629,MR1608365}. For other frameworks see
\cite{VermaSOIC} and references therein.


\subsection{General definitions}
\label{sec:defs}

Without loss of generality, a multiobjective optimization problem (MOP)
with $m$ objectives and $n$ decision variables can be formulated as follows:
\begin{equation}
\label{mo:problem}
\begin{array}{rl}
\text{minimize:} & \boldsymbol{f}(\boldsymbol{x})
=(f_1(\boldsymbol{x}),f_2(\boldsymbol{x}),\ldots,f_m(\boldsymbol{x}))^{\text{T}} \\
\text{subject to:} & \boldsymbol{x} \in \Omega
\end{array}
\end{equation}
where $\boldsymbol{x}$ is the decision vector, $\Omega \subseteq \mathbb{R}^n$
is the feasible decision space, and $\boldsymbol{f}(\boldsymbol{x})$ is the
objective vector defined in the attainable objective space $\Theta \subseteq \mathbb{R}^m$.

In multiobjective optimization, the Pareto dominance relation is usually used to define
the concepts of optimality. For two solutions $\boldsymbol{x}$ and $\boldsymbol{y}$
from $\Omega$, a solution $\boldsymbol{x}$ is said to dominate a solution $\boldsymbol{y}$
(denoted by $\boldsymbol{x} \prec \boldsymbol{y}$) if for all $i \in \{ 1, \ldots ,m\}$
$f_i (\boldsymbol{x}) \le f_i (\boldsymbol{y})$ and there exists at least one objective
such that $f_j (\boldsymbol{x}) < f_j (\boldsymbol{y})$.

A solution $\boldsymbol{x}^{\ast} \in \Omega$ is Pareto optimal if and only if:
\begin{equation*}
\nexists \boldsymbol{y} \in \Omega :\boldsymbol{y} \prec \boldsymbol{x}^{\ast}.
\end{equation*}

In the presence of multiple conflicting objectives, there is a set of optimal
solutions, known as the Pareto optimal set. For MOP~\eqref{mo:problem},
the Pareto optimal set (or Pareto set for short) is defined as:
\begin{equation*}
\mathcal{PS} = \{ \boldsymbol{x}^{\ast} \in \Omega \, | \, \nexists \boldsymbol{y}
\in \Omega : \boldsymbol{y} \prec \boldsymbol{x}^{\ast} \}.
\end{equation*}

For MOP~\eqref{mo:problem} and the Pareto set $\mathcal{PS}$, the Pareto optimal
front (or Pareto front for short) is defined as:
\begin{equation*}
\mathcal{PF} = \{ \boldsymbol{f}(\boldsymbol{x}^{\ast})
\in \Theta \, | \, \boldsymbol{x}^{\ast} \in \mathcal{PS} \}.
\end{equation*}
Since it is often not possible to obtain the whole Pareto set,
solving~\eqref{mo:problem} is usually understood as approximating
the Pareto set by obtaining a set of solutions that are as close
as possible to the Pareto set and as diverse as possible.


\subsection{Scalarization methods}
\label{sec:methods}

In the following, four scalarization methods for multiobjective optimization,
able to deal with convex and nonconvex Pareto fronts, are discussed.


\subsubsection{The $\epsilon$-Constraint method}
\label{sec:eps}

In the $\epsilon$-constraint method~\cite{Haimes1971}, one of the objective
functions is selected to be minimized, whereas all the other functions
are converted into constrains by setting an upper bound to each of them.
It can be defined as:
\begin{equation*}
\begin{array}{rlll}
\underset{\boldsymbol{x} \in \Omega}{\text{minimize:}}
& f_l(\boldsymbol{x}) & \\
\text{subject to:}   & f_i(\boldsymbol{x})\leq \epsilon_i,
& \forall i \in \{1,\ldots,m\} \wedge i \neq l, \\
\end{array}
\end{equation*}
where the $l$th objective is minimized and the parameter
$\epsilon_i$ represents an upper bound of the value of $f_i$.


\subsubsection{Chebyshev's method}
\label{sec:chb}

Chebyshev's method~\cite{Bowman1976} minimizes the weighted metric $L_p$,
with $p=\infty$, which measures the distance from any solution
to some reference point. It can be formulated as:
\begin{equation}
\label{method:chb}
\underset{\boldsymbol{x} \in \Omega}{\text{minimize:}}
\max \limits_{1\leq i\leq m} \, \{ w_i (f_i(\boldsymbol{x})-z^{\ast}_i) \},
\end{equation}
where $\boldsymbol{z}^{\ast} = (z^{\ast}_1,\ldots,z^{\ast}_m)^{\text{T}}$
is a reference point and $\boldsymbol{w} = (w_1, \ldots, w_m )^{\text{T}}$
is a weight vector such that $\sum\limits_{j=1}\limits^{m} w_j = 1$.


\subsubsection{The goal attainment method}
\label{sec:gam}

The problem shown in~\eqref{method:chb} can be reformulated as~\cite{MiettinenBook}:
\begin{equation*}
\left\{
\begin{array}{rl}
\underset{\boldsymbol{x} \in \Omega, \, \alpha \geq 0}{\text{minimize:}} & \alpha \\
\text{subject to:} & w_1(f_1(\boldsymbol{x}) - z^{\ast}_1) \leq \alpha, \\
& \qquad \vdots \\
& w_m(f_m(\boldsymbol{x}) - z^{\ast}_m) \leq \alpha, \\
\end{array}
\right.
\end{equation*}
where $\boldsymbol{z}^{\ast} = (z^{\ast}_1,\ldots,z^{\ast}_m)^{\text{T}}$
is a reference point, $\boldsymbol{w} = (w_1, \ldots, w_m )^{\text{T}}$
is a weight vector ($\sum\limits_{j=1}\limits^{m} w_j = 1$) and
$\boldsymbol{x} \in \Omega, \alpha \in \mathbb{R}_{+}$ are variables.
This method is often referred to as the goal attainment method~\cite{MiettinenBook}
or the Pascoletti-Serafini scalarization~\cite{PaSe84}.


\subsubsection{The normal constraint method}
\label{sec:nc}

The normal constraint (NC) method~\cite{Messac2003,Messac2004} minimizes one
of the objective functions and uses an inequality constraint reduction
of the feasible objective space. It can be formulated as:
\begin{equation*}
\begin{array}{rll}
\underset{\boldsymbol{x} \in \Omega}{\text{minimize:}}
& \overline{f}_l(\boldsymbol{x}) & \\
\text{subject to:}
& \boldsymbol{\overline{v}}_i^{\text{T}}(\boldsymbol{\overline{f}}(\boldsymbol{x})
-\boldsymbol{z})\leq 0, & \forall i \in \{1,\ldots, m \} \wedge i \neq l,
\end{array}
\end{equation*}
where $\boldsymbol{\overline{v}}_i$ is a vector from the $i$th corner,
$\boldsymbol{\overline{\mu}}^{i\ast}$, to the fixed $l$th corner,
$\boldsymbol{\overline{\mu}}^{l\ast}$, of the Pareto front
\[
\boldsymbol{\overline{v}}_i=\boldsymbol{\overline{\mu}}^{l\ast}
-\boldsymbol{\overline{\mu}}^{i\ast},
\quad \forall i \in \{1,\ldots, m \} \wedge i \neq l;
\]
$\boldsymbol{z}$ is a point on the hyperplane corresponding
to a given vector $\boldsymbol{w} = (w_1, \ldots, w_m )^{\text{T}}$,
\[
\boldsymbol{z}= \boldsymbol{\Phi} \boldsymbol{w},
\]
where $\boldsymbol{\Phi}=(\overline{\boldsymbol{\mu}}^{1\ast},
\ldots,\overline{\boldsymbol{\mu}}^{m\ast})$ is a $m\times m$
matrix and $\sum\limits_{j=1}\limits^{m} w_j = 1$.

To cope with differently scaled objectives,
the NC method normalizes the objective vector as:
\begin{equation}
\label{normalization}
\overline{f}_i=\frac{f_i-z_i^{\text{ideal}}}{z_i^{\text{nadir}}
-z_i^{\text{ideal}}}, \quad \forall i \in \{1,\ldots, m \},
\end{equation}
where $z_i^{\text{ideal}}$ and $z_i^{\text{nadir}}$ are the $i$th
components of the ideal and nadir points, respectively.


\subsection{An algorithm for the Pareto set approximation}
\label{sec:solver}

The above discussed methods convert a multiobjective optimization
problem into a single-objective problem depending on some parameters.
Solving the corresponding problem for different parameter settings
allows to approximate multiple Pareto optimal solutions. The outline
of an approach to approximate the Pareto set of problem~\eqref{dengue:mop}
is shown in Algorithm~\ref{mo:solver}.

\begin{algorithm}
\caption{Pareto set approximation}\label{mo:solver}
\begin{algorithmic}[1]
\State initialize: $\boldsymbol{x}^{(0)}$, $B$
\State $A \gets \{ \}$;
\For{$\beta \in B$}
\State for initial point $\boldsymbol{x}^{(0)}$,
find $\boldsymbol{x}^{\ast}$ that minimizes $f_{\beta}(\boldsymbol{x}, \beta)$;
\State $A\gets A \cup \{\boldsymbol{x}^{\ast}\}$;
\State $\boldsymbol{x}^{(0)}\gets \boldsymbol{x}^{\ast}$;
\EndFor
\State output: $A$;
\end{algorithmic}
\end{algorithm}

In line~1 of Algorithm~\ref{mo:solver}, an initial point is generated
as the null vector, $\boldsymbol{x}^{(0)}=\boldsymbol{0}$, and a set
of parameters, $B$, used for scalarization, is initialized. For the
$\epsilon$-constraint method, $\beta \in B$ corresponds to the value
of $\epsilon$. In the case of the Chebyshev and goal attainment methods,
$\beta$ is a tuple of the form $\{ \boldsymbol{w}, \boldsymbol{z}^{\ast}\}$,
whereas $\beta$ is a tuple of the form
$\{ \boldsymbol{w}, \boldsymbol{z}^{\text{ideal}},
\boldsymbol{z}^{\text{nadir}}\}$ for the normal constraint method. In lines 3--7,
for each corresponding scalarizing function, $f_{\beta}(\boldsymbol{x}, \beta)$,
a minimizer, $\boldsymbol{x}^{\ast}$, is found using a single-objective
optimizer and added to the approximation set, $A$. The minimizer of the previous
scalarizing function becomes an initial point for solving the subsequent problem,
to ensure efficiency of the approach.


\section{Numerical experiments}
\label{sec:experiments}

This section provides results of the numerical experiments conducted
for the dengue transmission model. The comparison of different
scalarization methods is performed. The analysis of obtained
solutions and the dengue dynamics are presented.


\subsection{Experimental setup}
\label{sec:setup}

The fourth-order Runge--Kutta method is used to discretize the control
and state variables of the system~\eqref{dengue:model}. The period
$[0,84]$ of 84 days is discretized using the equally spaced time intervals of
$0.25$ (6 hours). This results in an optimization problem with $\boldsymbol{x}
\in [0,1]^{337}$, where $\boldsymbol{x}$ denotes the discrete control
in~\eqref{dengue:model}. The integrals defining the objective functionals
in~\eqref{dengue:mop} are calculated using the trapezoidal rule.

For approximating the Pareto set of~\eqref{dengue:mop}, scalarization
is performed by defining and solving 100 scalarizing problems corresponding
to the above discussed methods. As a single-objective optimizer
(line~4 in Algorithm~\ref{mo:solver}), the MATLAB\textsuperscript{\textregistered}
function \texttt{fmincon} with a sequential quadratic programming algorithm is used,
setting the maximum number of function evaluations to $20,000$.


\subsection{Different scalarization methods}
\label{sec:comparison}

For comparison, the outcomes obtained by different scalarization
methods are assessed using the hypervolume~\cite{ZitzlerThiele1998}.
The hypervolume can be defined as the Lebesgue measure, $\Lambda$,
of the union of hypercuboids in the objective space:
\begin{equation*}
HV=\Lambda\left(\bigcup_{\boldsymbol{a} \in A \wedge \boldsymbol{r}}
[ f_1(\boldsymbol{a}),r_1]\times \cdots \times [f_m(\boldsymbol{a}),r_m] \right),
\end{equation*}
\makebox[\linewidth][s]{where $A=\{ \boldsymbol{a}_1,\ldots,\boldsymbol{a}_{|A|} \}$
denotes a set of nondominated solutions and} \\
$\boldsymbol{r}=(r_1,\ldots,r_m)^{\text{T}}$ is a reference point.
The hypervolume calculates a portion of the objective space dominated by $A$.
It can measure both convergence and diversity. The higher the value of $HV$,
the better the quality of $A$. For calculating the hypervolume, the objectives
are normalized using~\eqref{normalization} and $[1, 1]$ is used as a reference point.
\begin{table}
\centering
\small
\begin{tabular}{ccccc}
\toprule
\backslashbox{measure}{method}& $\epsilon$-Constraint
& Chebyshev & Goal attainment & Normal constraint \\ \hline
Hypervolume & 0.972798 & 0.959887 & 0.959889 & {\bf 0.977248} \\
\bottomrule	
\end{tabular}
\caption{Hypervolume values for trade-off solutions obtained using different
scalarization methods (the higher the better). The best value is marked bold.}
\label{tab:hv}
\end{table}
\begin{figure}
\centering
\subfloat[$\epsilon$-Constraint method]{\epsfig{file=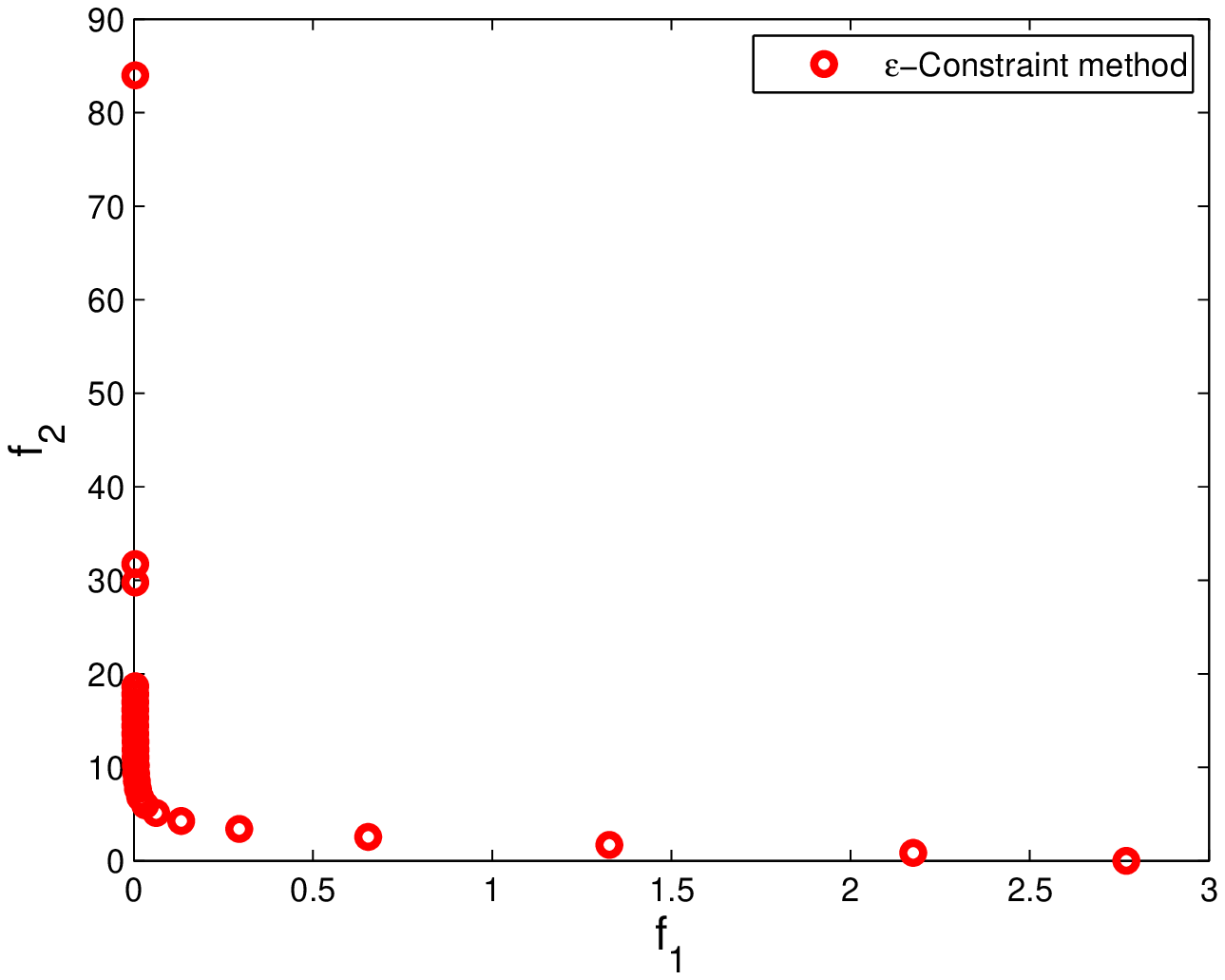,width=0.45\textwidth}\label{fig:methods:eps}}%
\subfloat[Chebyshev method]{\epsfig{file=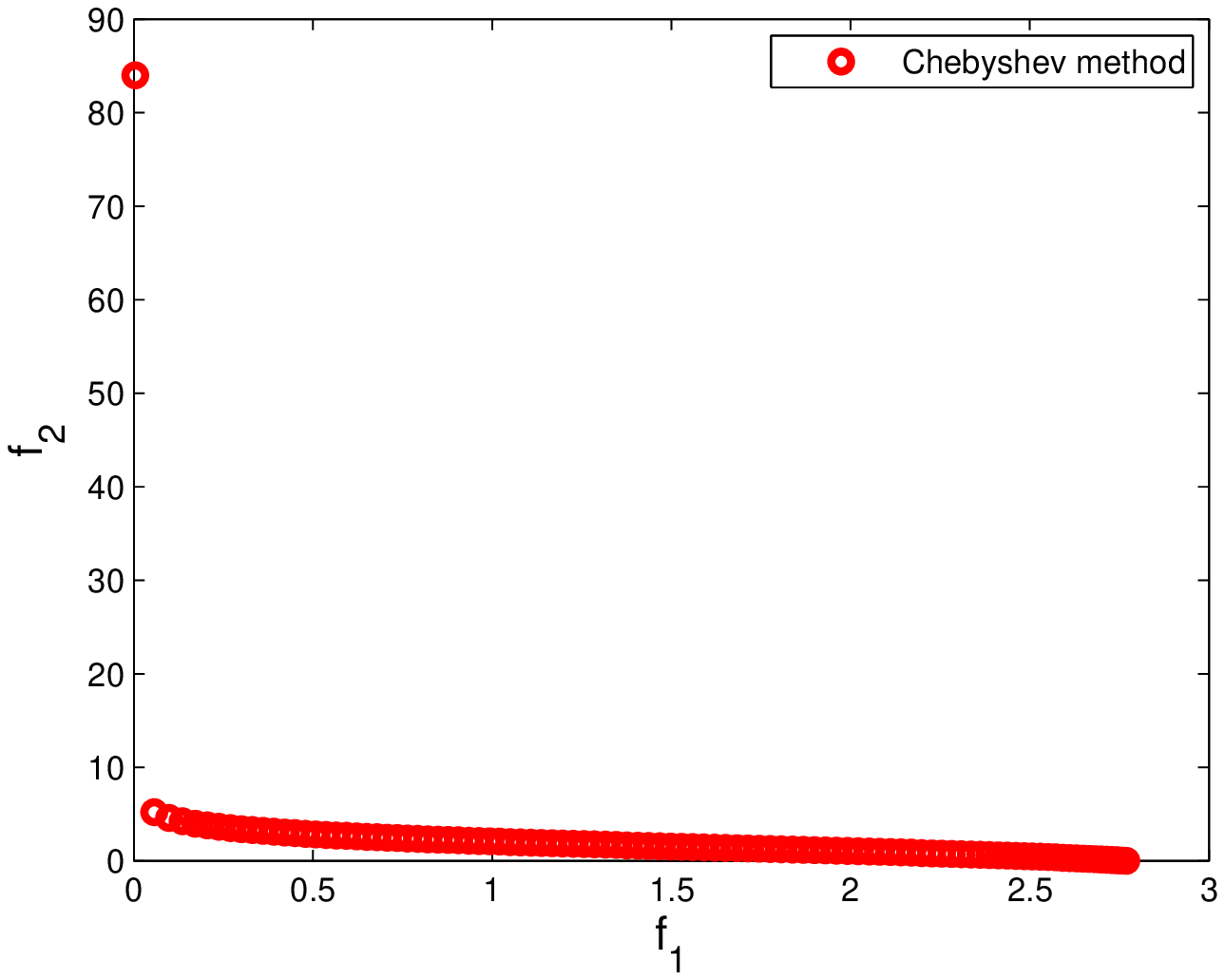,width=0.45\textwidth}\label{fig:methods:chb}}\\
\subfloat[Goal attainment method]{\epsfig{file=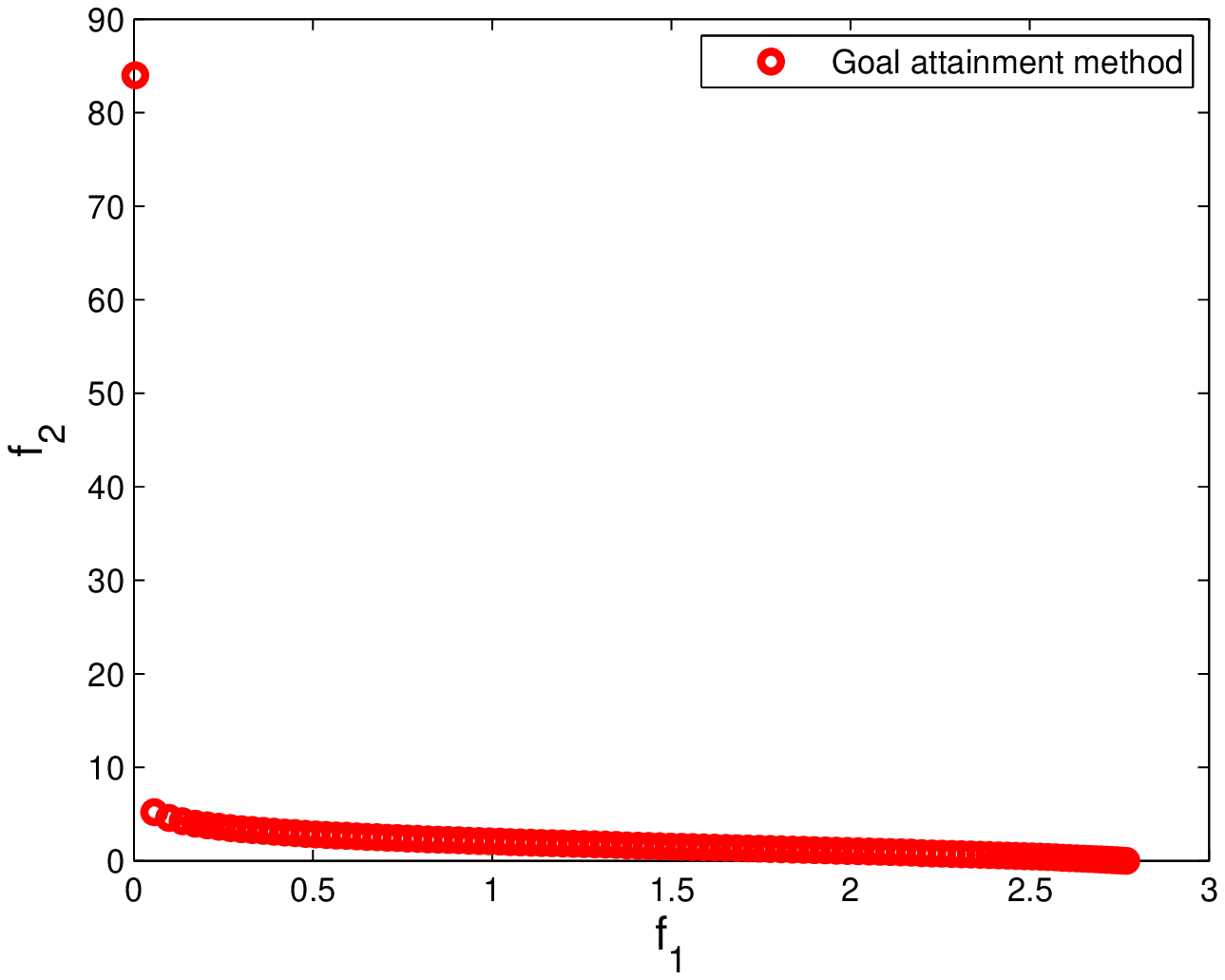,width=0.45\textwidth}\label{fig:methods:gam}}%
\subfloat[Normal constraint method]{\epsfig{file=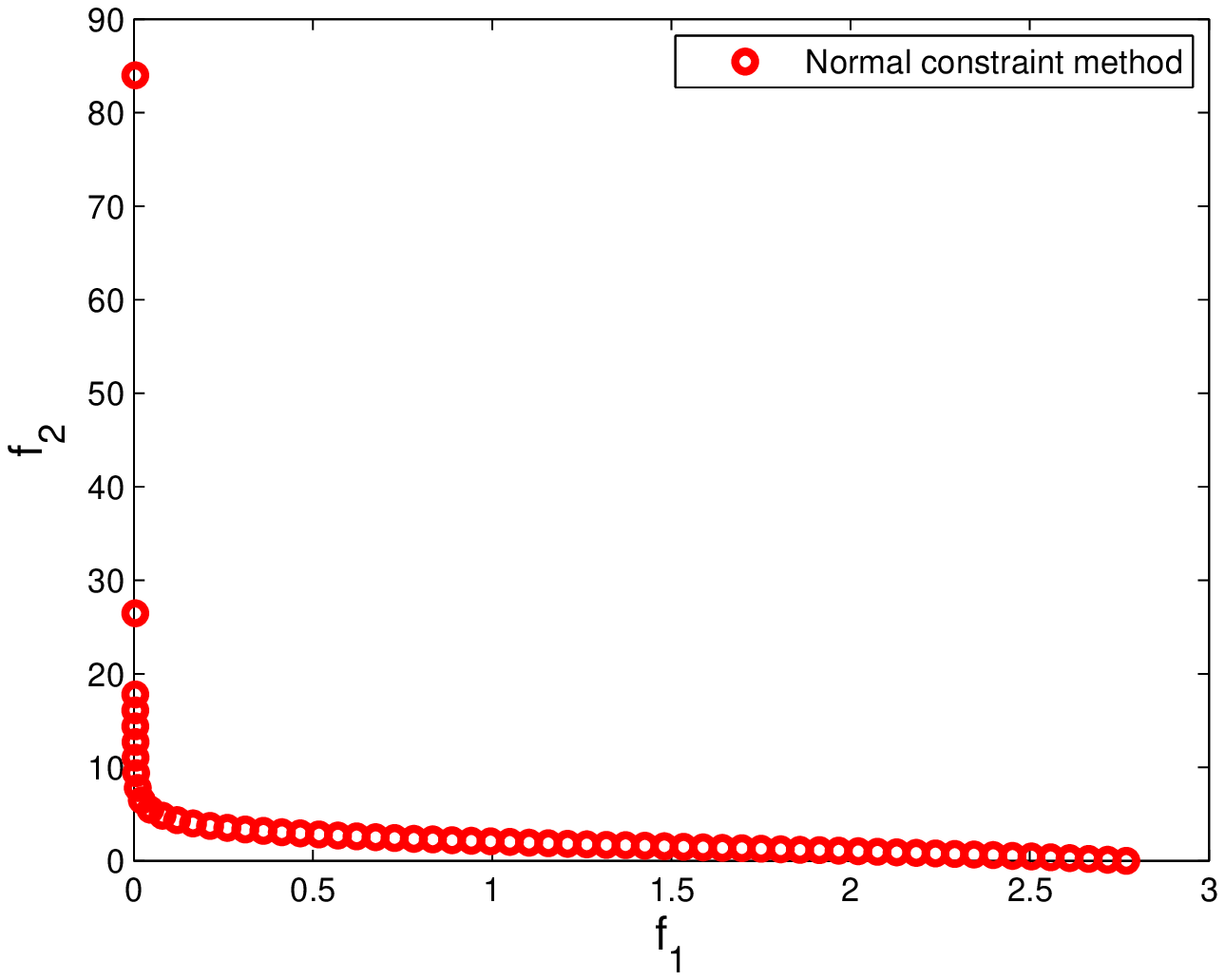,width=0.45\textwidth}\label{fig:methods:nc}}\\
\caption{Trade-off curves obtained by different methods.}\label{fig:methods}
\end{figure}

Table~\ref{tab:hv} presents the results for the four scalarization methods
with respect to the hypervolume. The best result is obtained by the normal
constraint method followed by the $\epsilon$-constraint method. Both methods
minimize $f_1$, when $f_2$ is used to define the constraint for the objective
space reduction. Since the NC method uses evenly distributed points on the
ideal plane, the final points on the Pareto front approximation are likely
to be more evenly distributed than using the usual $\epsilon$-constraint method.
This allows to achieve the highest hypervolume. The Chebyshev method provides
a slightly worse result than the goal attainment method. This can be because
the scalarizing function used by the Chebyshev method is nondifferentiable,
which usually poses additional difficulties for optimization.

Figure~\ref{fig:methods} shows the Pareto front approximations obtained by the methods.
The plots confirm the previous observations based on the hypervolume values.
The Chebyshev and goal attainment methods provide visually similar results.
The $\epsilon$-constraint method attempts to find evenly distributed points
with respect to the uniform division of $f_2$, whereas the NC method seeks
a uniform distribution of points, according to points on the hyperplane passing
through the corner points of the Pareto front. Since the NC method gives clearly
better results for~\eqref{dengue:mop}, solutions obtained by this method
will be employed for the analysis of the dengue transmission in what follows.


\subsection{Dengue dynamics}
\label{sec:dynamics}

\begin{figure}
\centering
\subfloat[$0\leq f_2 \leq 10$]{\epsfig{file=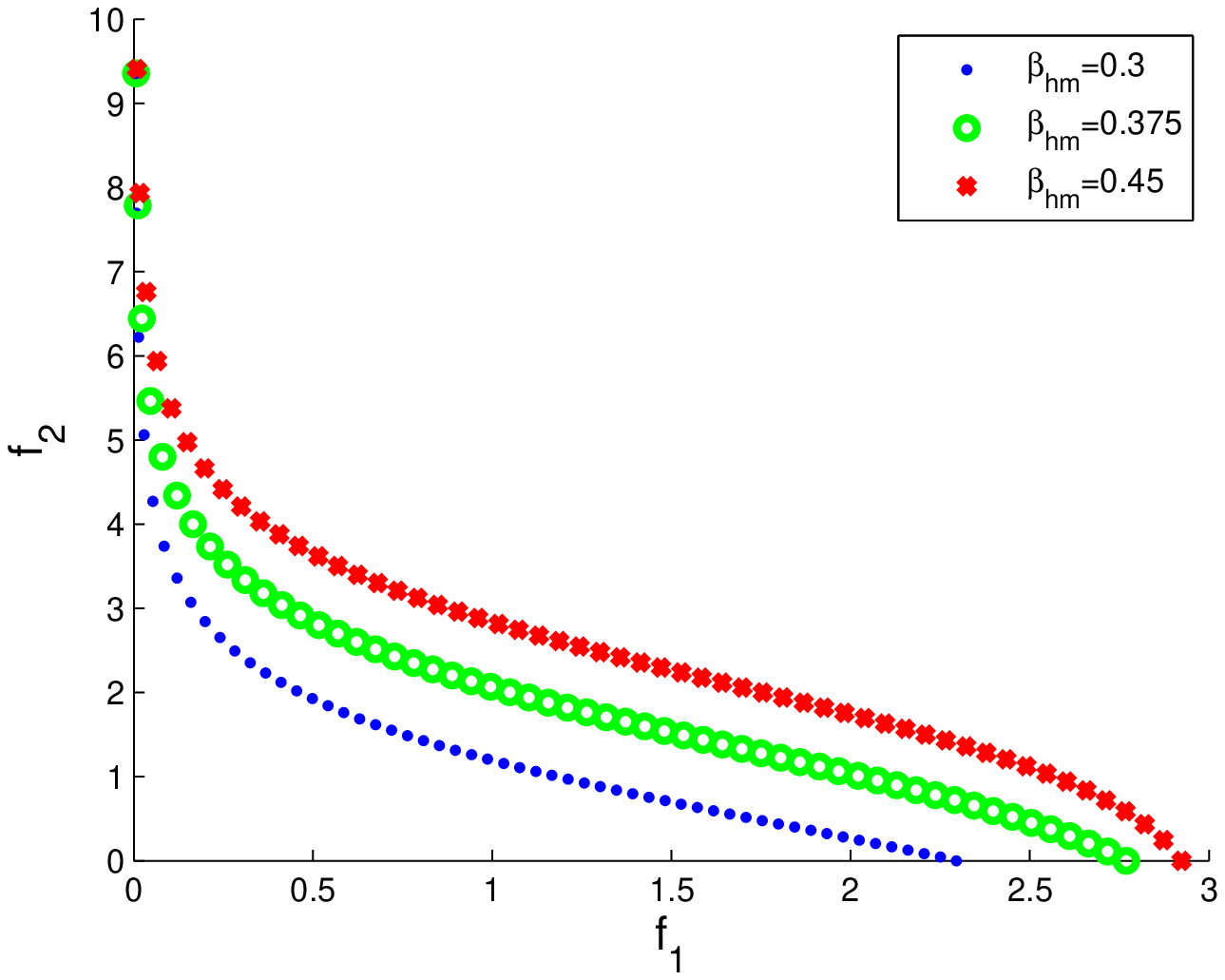,width=0.45\textwidth}\label{fig:betahms:1}}%
\subfloat[$10\leq f_2 \leq 90$]{\epsfig{file=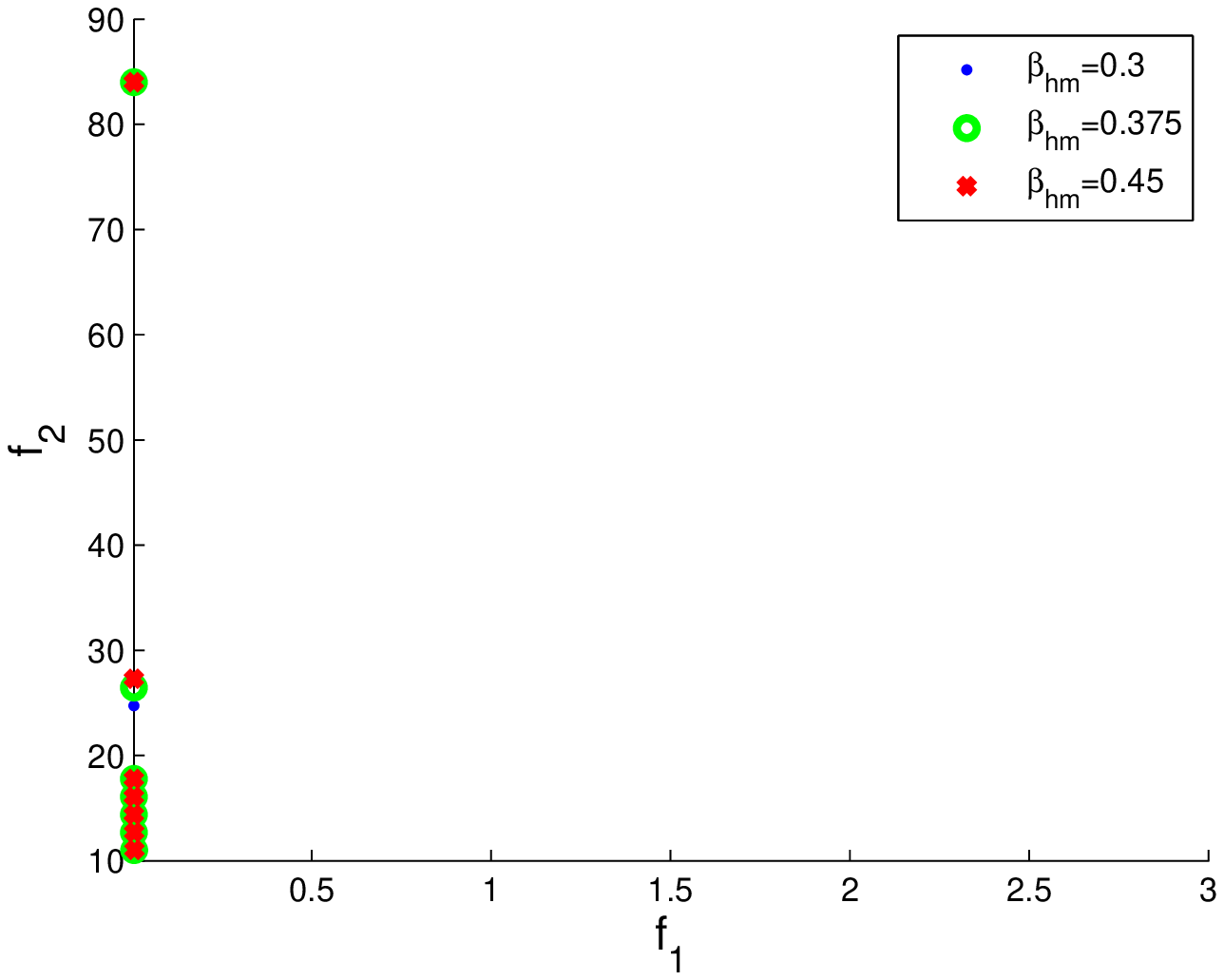,width=0.45\textwidth}\label{fig:betahms:2}}\\
\caption{Trade-off curves for different values of $\beta_{hm}$.}\label{fig:betahms}
\end{figure}

\begin{figure}
\centering
\subfloat[$0\leq f_2 \leq 10$]{\epsfig{file=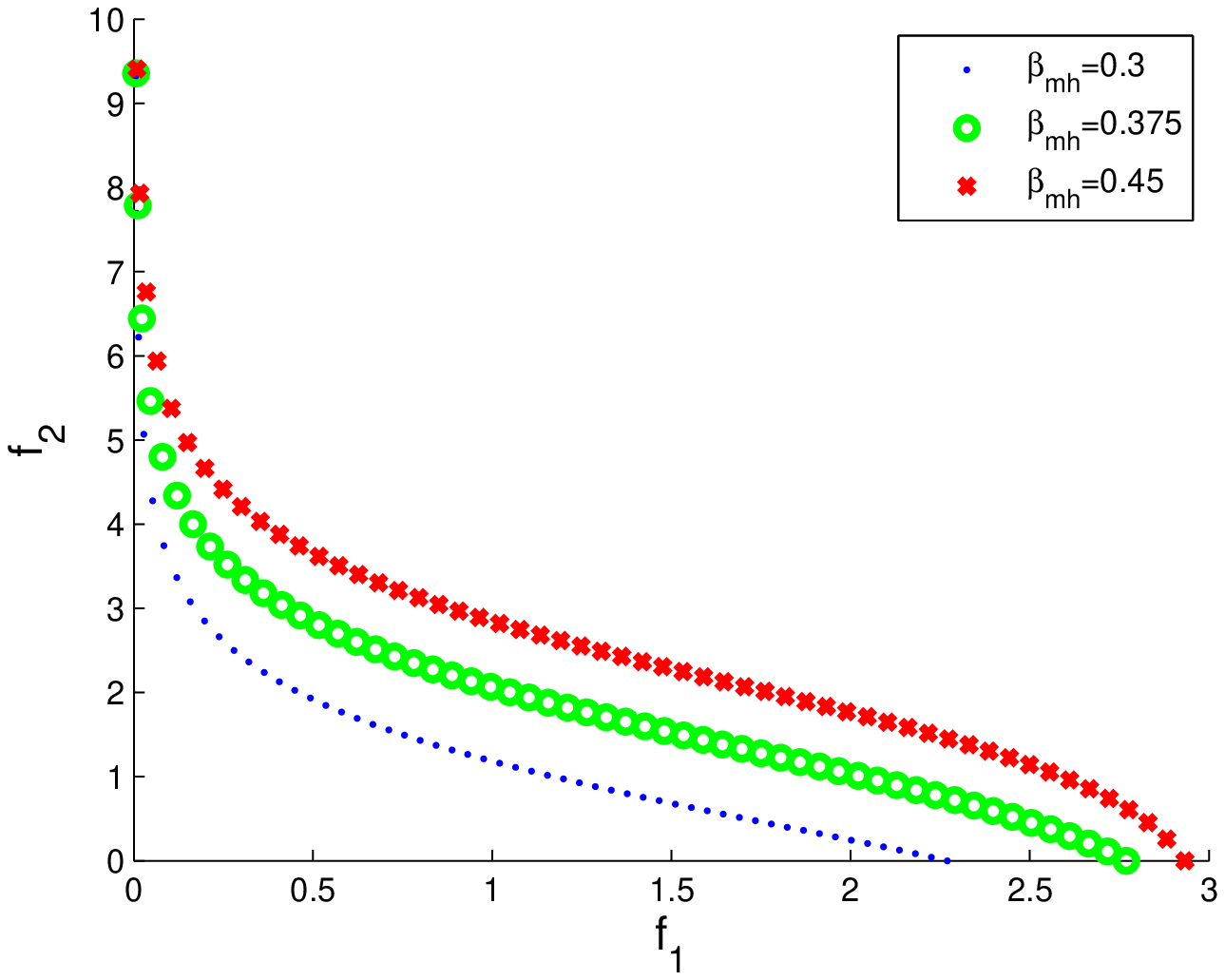,width=0.45\textwidth}\label{fig:betamhs:1}}%
\subfloat[$10\leq f_2 \leq 90$]{\epsfig{file=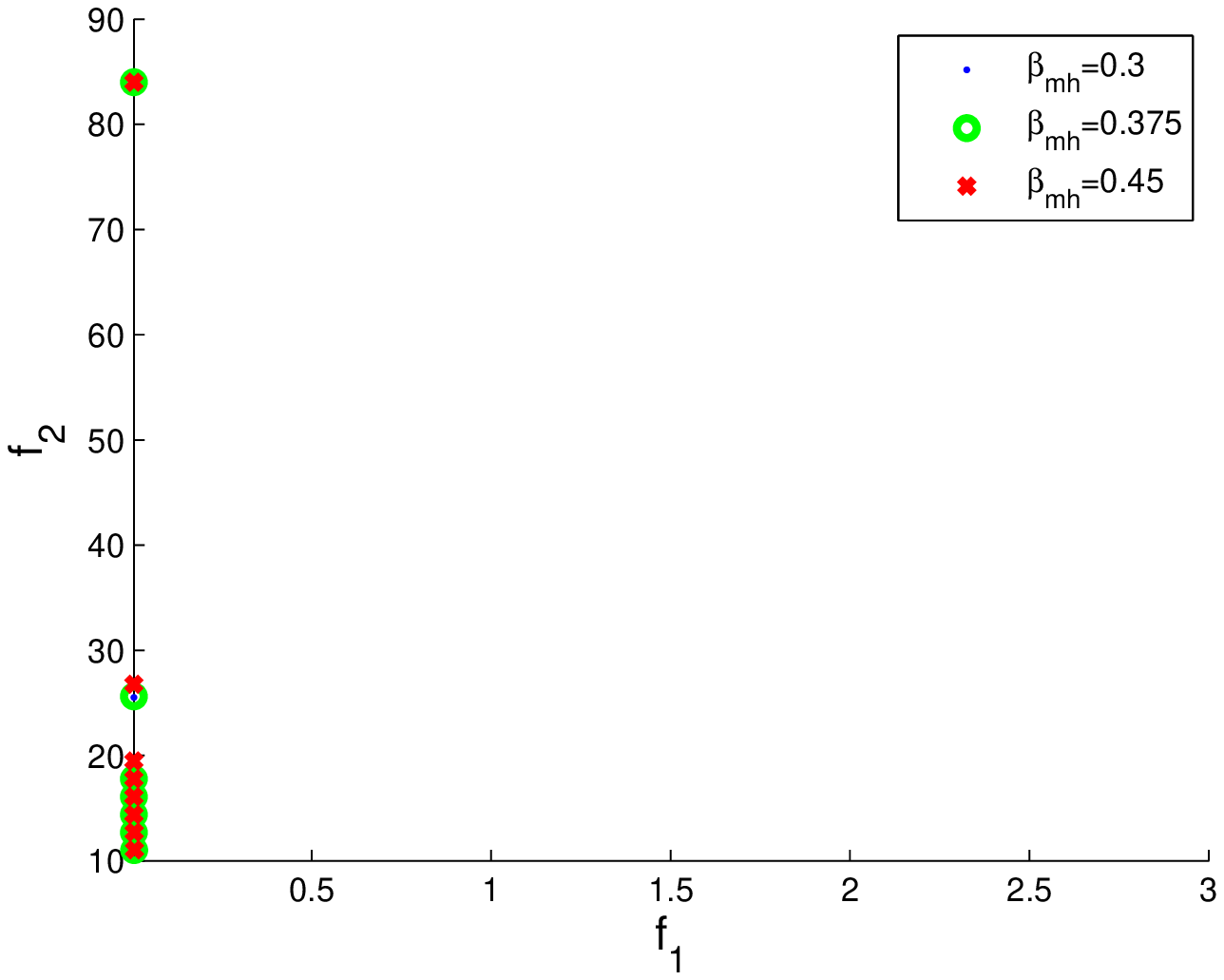,width=0.45\textwidth}\label{fig:betamhs:2}}\\
\caption{Trade-off curves for different values of $\beta_{mh}$.}\label{fig:betamhs}
\end{figure}

Figures~\ref{fig:betahms} and~\ref{fig:betamhs} present the trade-off solutions
obtained for different values of transition probabilities $\beta_{hm}$ and $\beta_{mh}$,
respectively. As it was seen during discussion of different scalarization methods,
the trade-off curve for~\eqref{dengue:mop} exhibits a smooth relation between
the insecticide cost, $f_2$, and the infected human population, $f_1$, in the
range $0\leq f_2\leq 5$. This suggests that implementing optimal control strategies
in this range can produce the most significant reduction in infected individuals,
having the highest ratio between the desirable effect and the cost. However,
starting from some further point, reducing the number of infected humans can be
achieved through exponential increase in spendings for insecticide. Even a small
decrease corresponds to a high increase in expenses for insecticide. Scenarios
represented by this part of the trade-off curve can be unacceptable from the
economical point of view. Furthermore, it should be noted that even with the
maximum spending it is not possible to eradicate the disease, with the lowest
obtained value of the infected population being 0.0042.

\begin{figure}
\centering
\subfloat[$i_h$ for $c(\cdot) \equiv 0$ and varying $\beta_{hm}$]{%
\epsfig{file=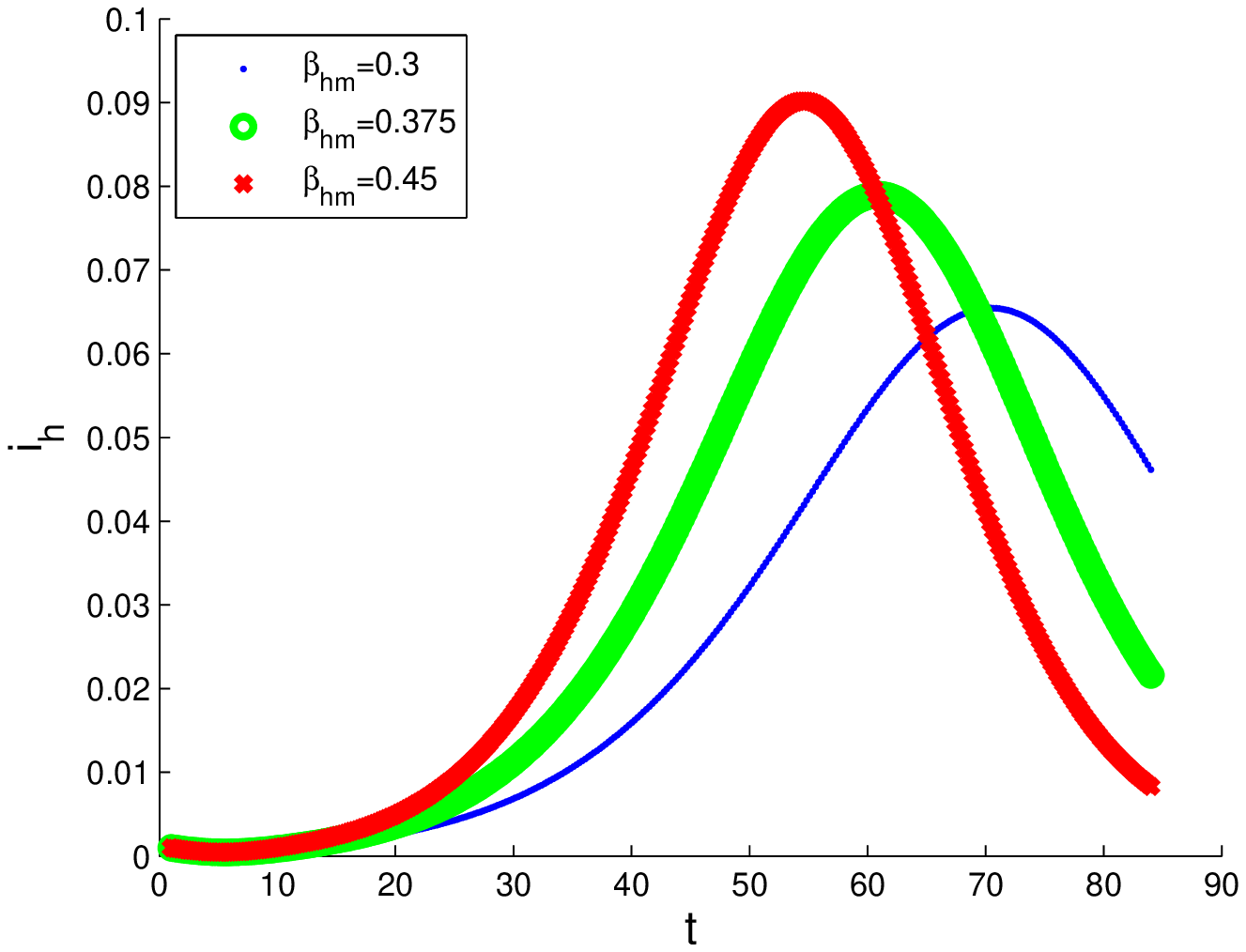,width=0.45\textwidth}\label{fig:dynamics:1}}%
\subfloat[$i_h$ for $c(\cdot) \equiv 0$ and varying $\beta_{mh}$]{%
\epsfig{file=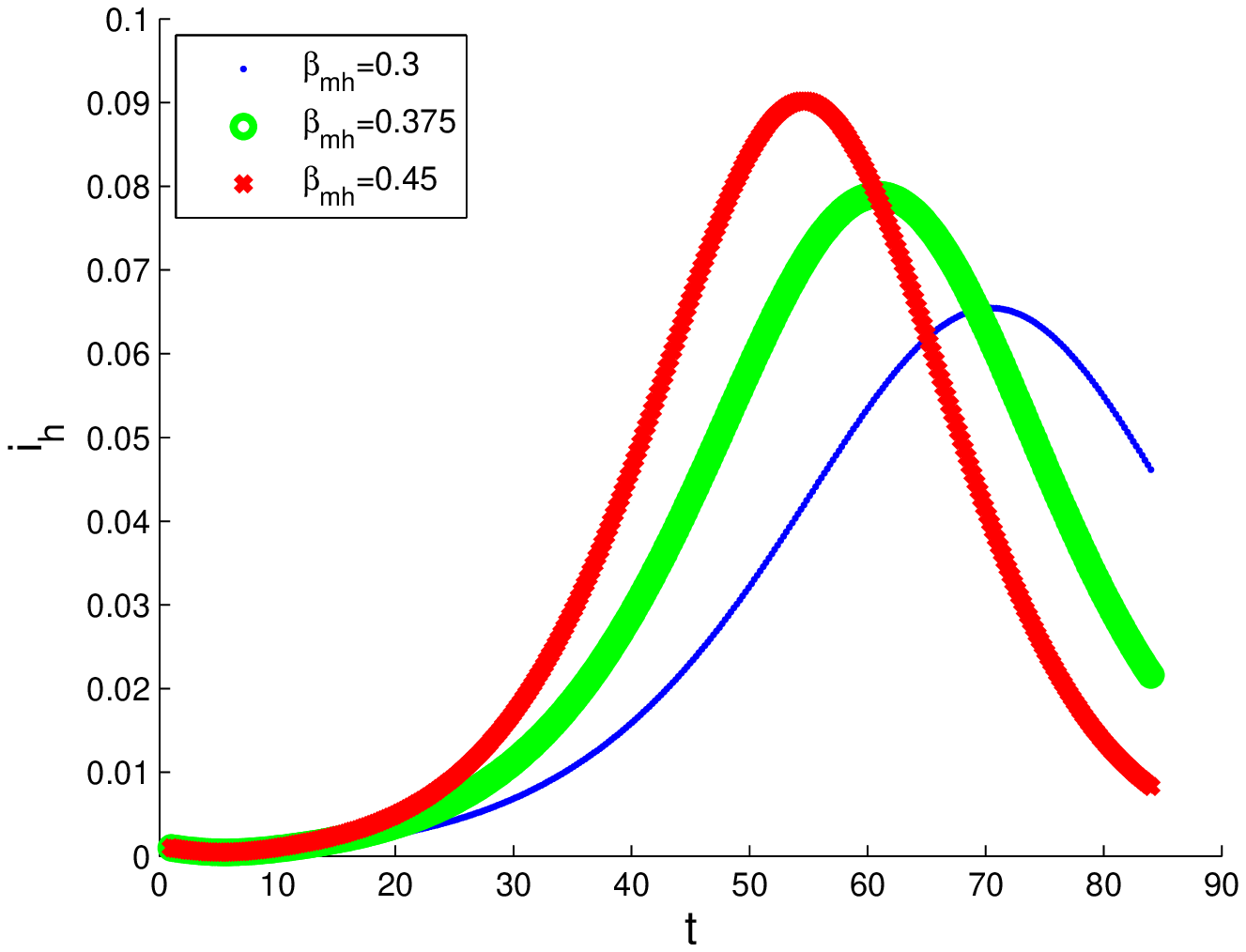,width=0.45\textwidth}\label{fig:dynamics:2}}\\
\subfloat[$i_h$ for $c(\cdot) \equiv 1$ and varying $\beta_{hm}$]{%
\epsfig{file=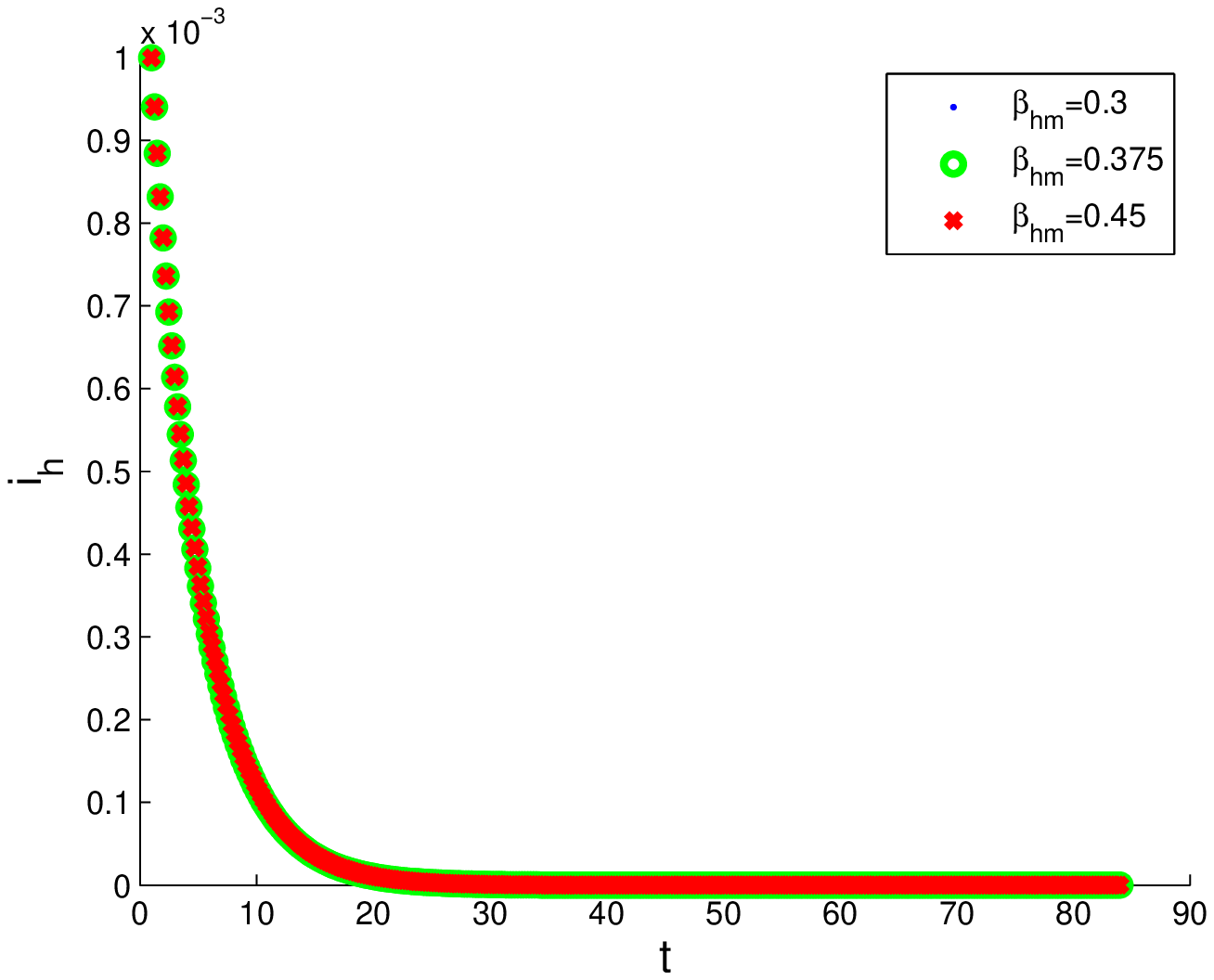,width=0.45\textwidth}\label{fig:dynamics:3}}%
\subfloat[$i_h$ for $c(\cdot) \equiv 1$ and varying $\beta_{mh}$]{%
\epsfig{file=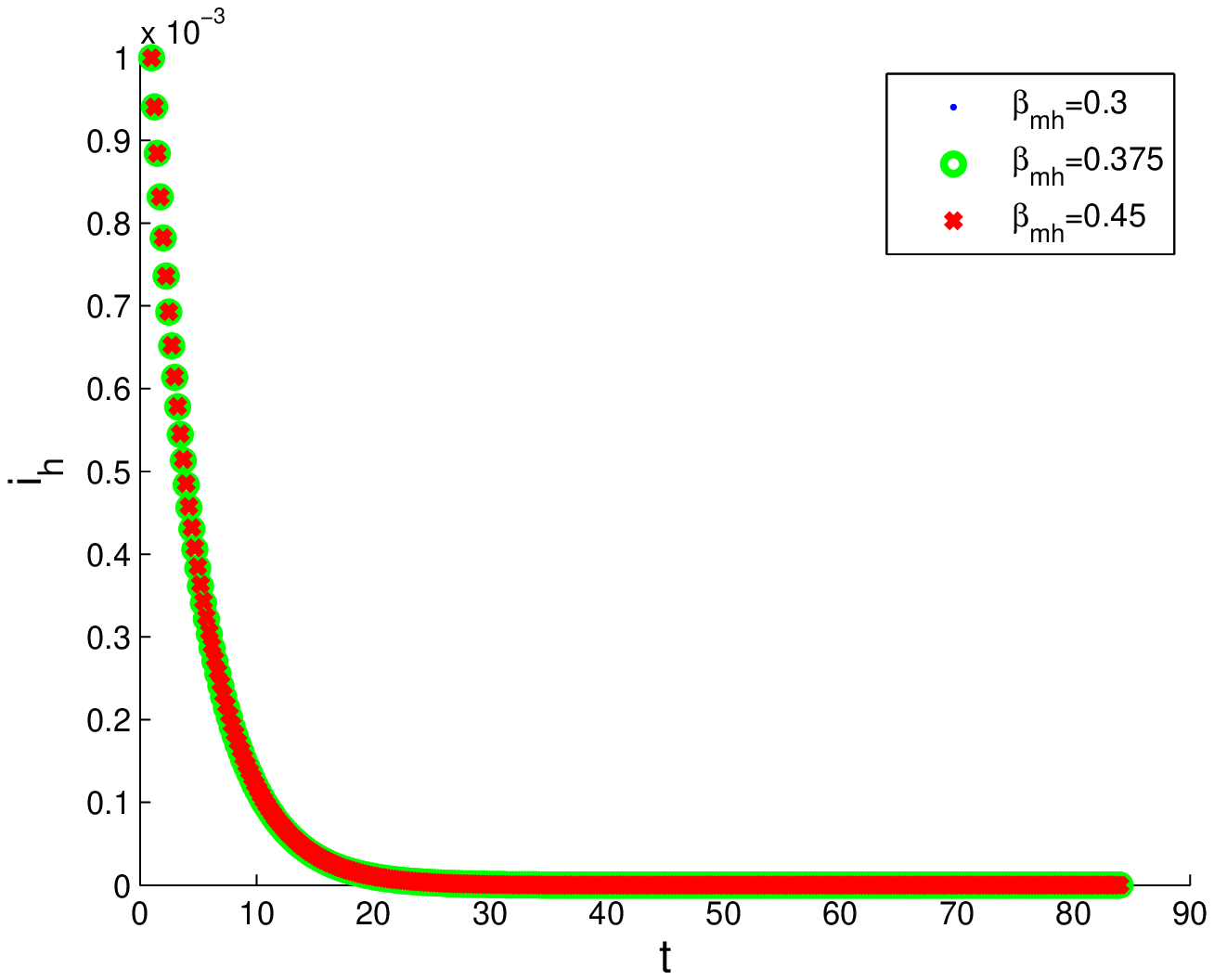,width=0.45\textwidth}\label{fig:dynamics:4}}\\
\caption{Dynamics of the state variable $i_h$ for varying values of $\beta_{hm}$
and $\beta_{mh}$, corresponding to extreme solutions.}\label{fig:dynamics}
\end{figure}

To provide a better visualization of the most interesting parts of the trade-off curves
in Figures~\ref{fig:betahms} and~\ref{fig:betamhs}, the plots are separated into two
parts and presented in different figures, according to the range of $f_2$.
From Figures~\ref{fig:betahms:1} and~\ref{fig:betamhs:1}, it can be seen that
the higher the corresponding transition probability, the larger number of infected
individuals. Though this difference reduces as more control is applied, vanishing
for the extreme scenario with $c(\cdot) \equiv 1$. Varying $\beta_{hm}$ and $\beta_{mh}$
produces apparently similar effects on the shape of the Pareto front. From
Figures~\ref{fig:betahms:1} and~\ref{fig:betamhs:1}, it can be observed that
for higher values of $\beta_{hm}$ and $\beta_{mh}$ the part of the trade-off curve
in $0\leq f_2\leq 5$ becomes increasingly nonconvex. Solutions in this region can be
unattainable for methods that face difficulties in dealing
with nonconvexities in the Pareto front.

Figure~\ref{fig:dynamics} illustrates dengue epidemics, corresponding to extreme scenarios
($c(\cdot) \equiv 0$ and $c(\cdot) \equiv 1$). It can be seen that for higher values of transition
probabilities the peak in infected population occurs early, having larger values
(Figures~\ref{fig:dynamics:1} and~\ref{fig:dynamics:2}).  On the other hand,
applying the maximum control, $c(\cdot) \equiv 1$, the dynamics of $i_h(t)$ remain
unchanged for different values of $\beta_{hm}$ and $\beta_{mh}$
(Figures~\ref{fig:dynamics:3} and~\ref{fig:dynamics:4}).

\begin{figure}
\centering
\includegraphics[width=0.45\textwidth]{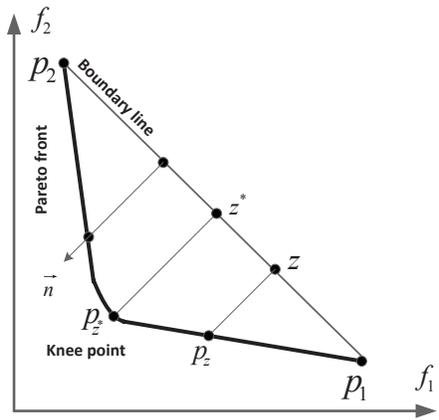}
\caption{Knee point identification.}\label{knee}
\end{figure}

To discuss intermediate scenarios, the concept of the knee solution of the Pareto
front is adopted. According to~\cite{Das99}, it can be defined as follows.
Given a Pareto front with the normalized objectives, a boundary line $L(p_1, p_2)$
is constructed though two extreme points $p_1$ and $p_2$. For any point $z$
on the boundary line $L(p_1, p_2)$, a point $p_z$ on the Pareto front along
the normal $\vec{n}$ of the boundary line is identified. The Pareto optimal
point $p_{z^{\ast}}$ with the maximum distance from its corresponding boundary
point $z^{\ast}$ along the normal direction is defined as the knee point.
This concept is illustrated in Figure~\ref{knee}. The knee solution of the
Pareto front is attractive, since it is often considered
as the optima in objective trade-offs.

\begin{figure}
\centering
\subfloat[$c$ for knee point and varying $\beta_{hm}$]{%
\epsfig{file=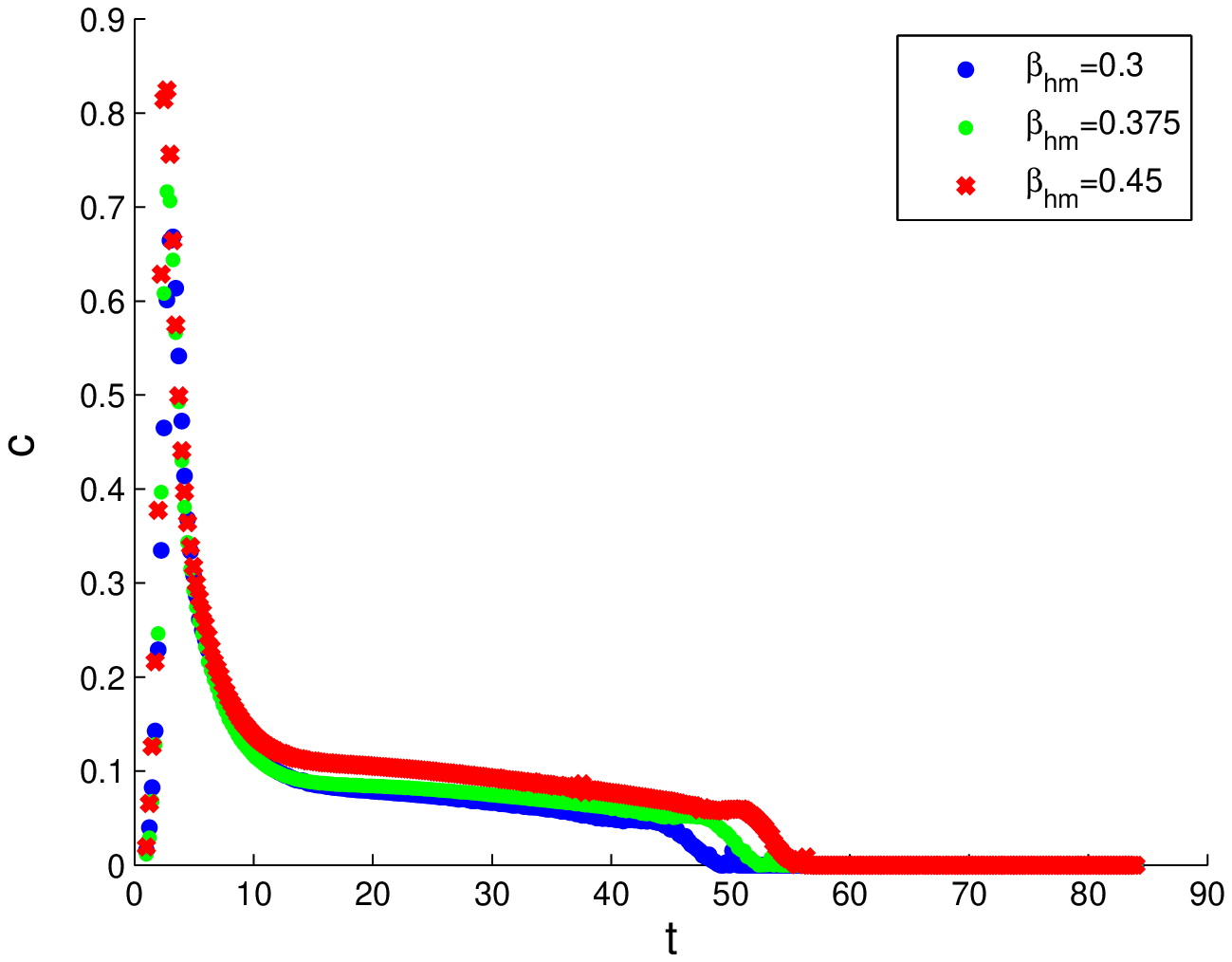,width=0.45\textwidth}\label{fig:dynamics:knee:1}}%
\subfloat[$c$ for knee point and varying $\beta_{mh}$]{%
\epsfig{file=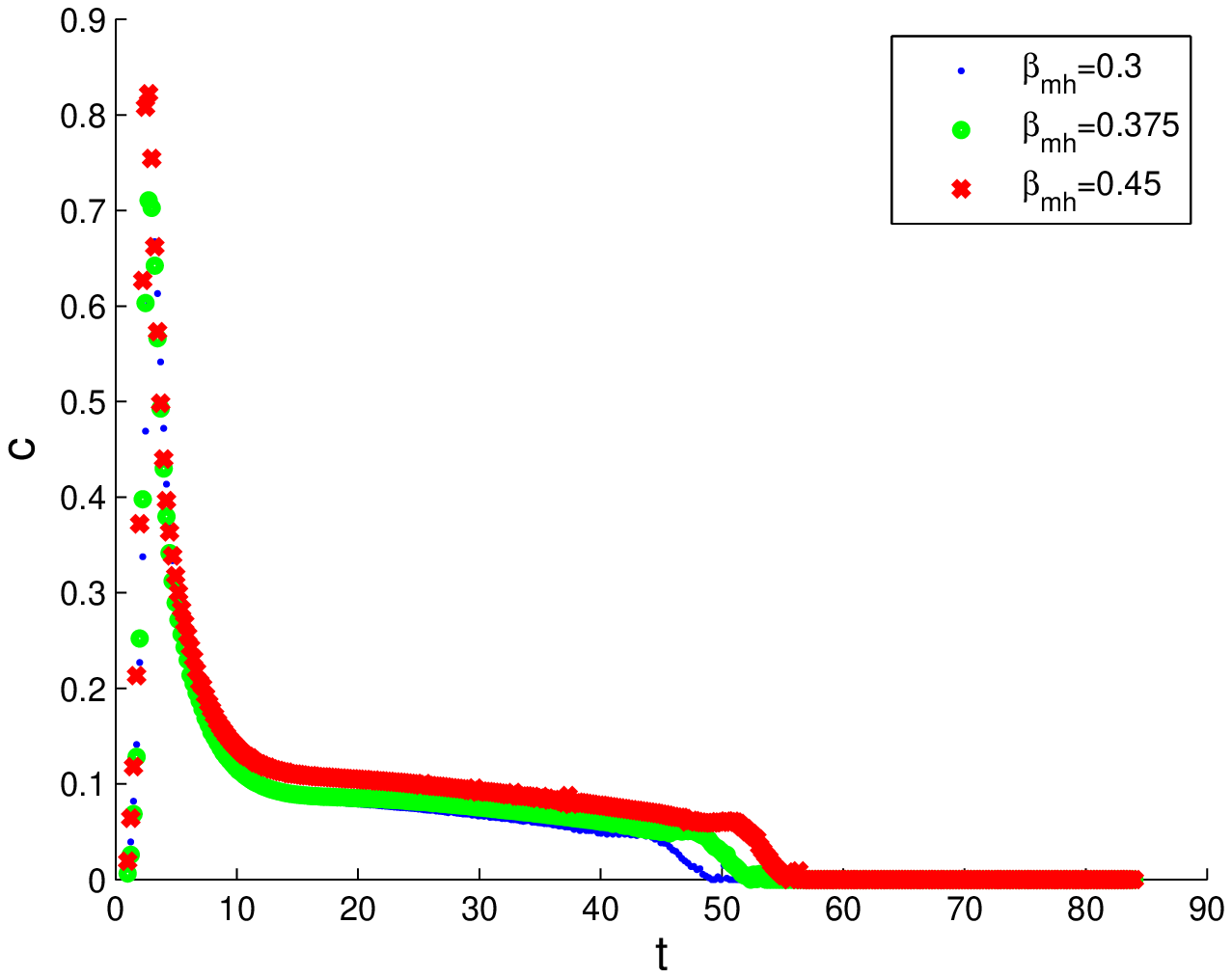,width=0.45\textwidth}\label{fig:dynamics:knee:2}}\\
\subfloat[$i_h$ for knee point and varying $\beta_{hm}$]{%
\epsfig{file=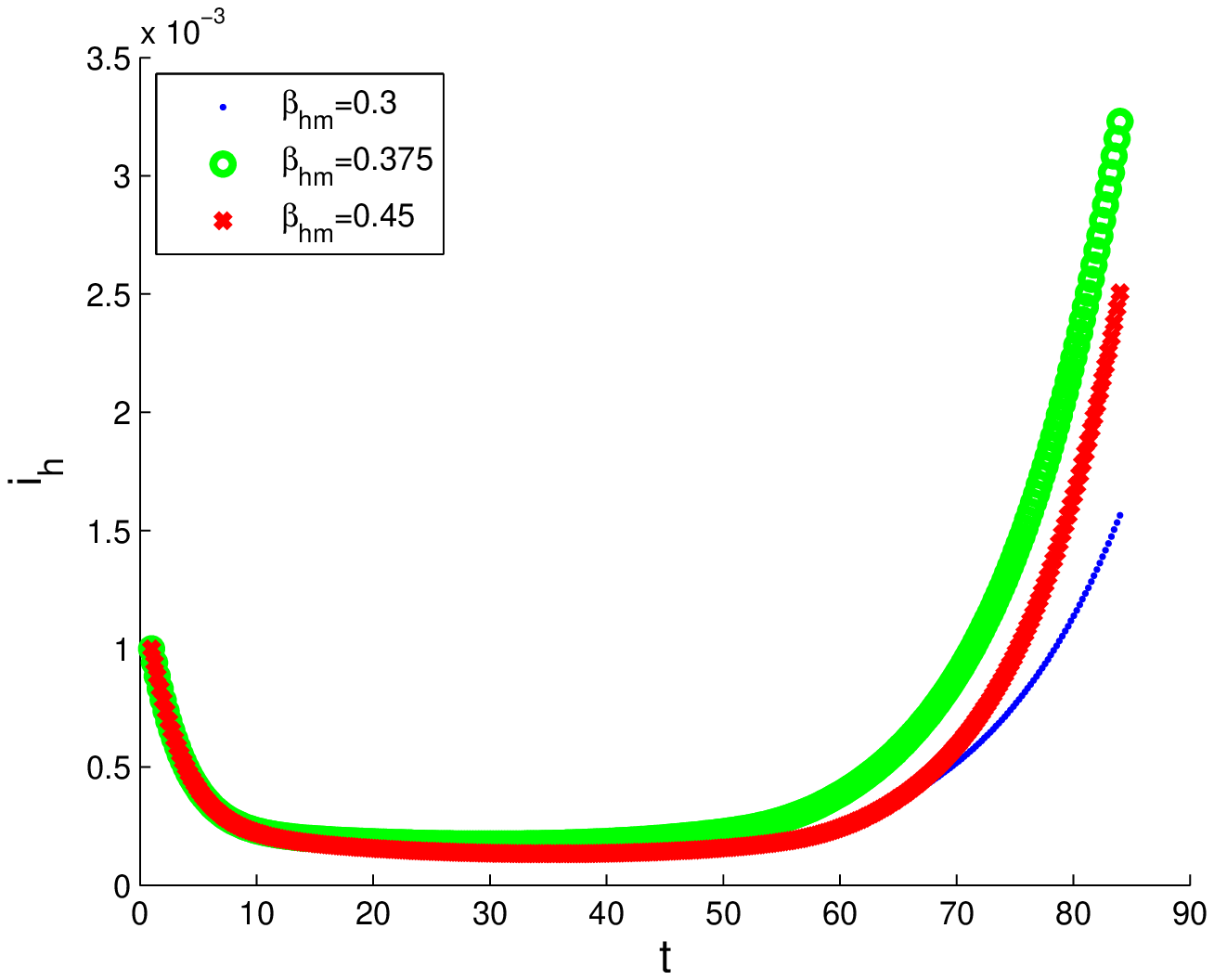,width=0.45\textwidth}\label{fig:dynamics:knee:3}}%
\subfloat[$i_h$ for knee point and varying $\beta_{mh}$]{%
\epsfig{file=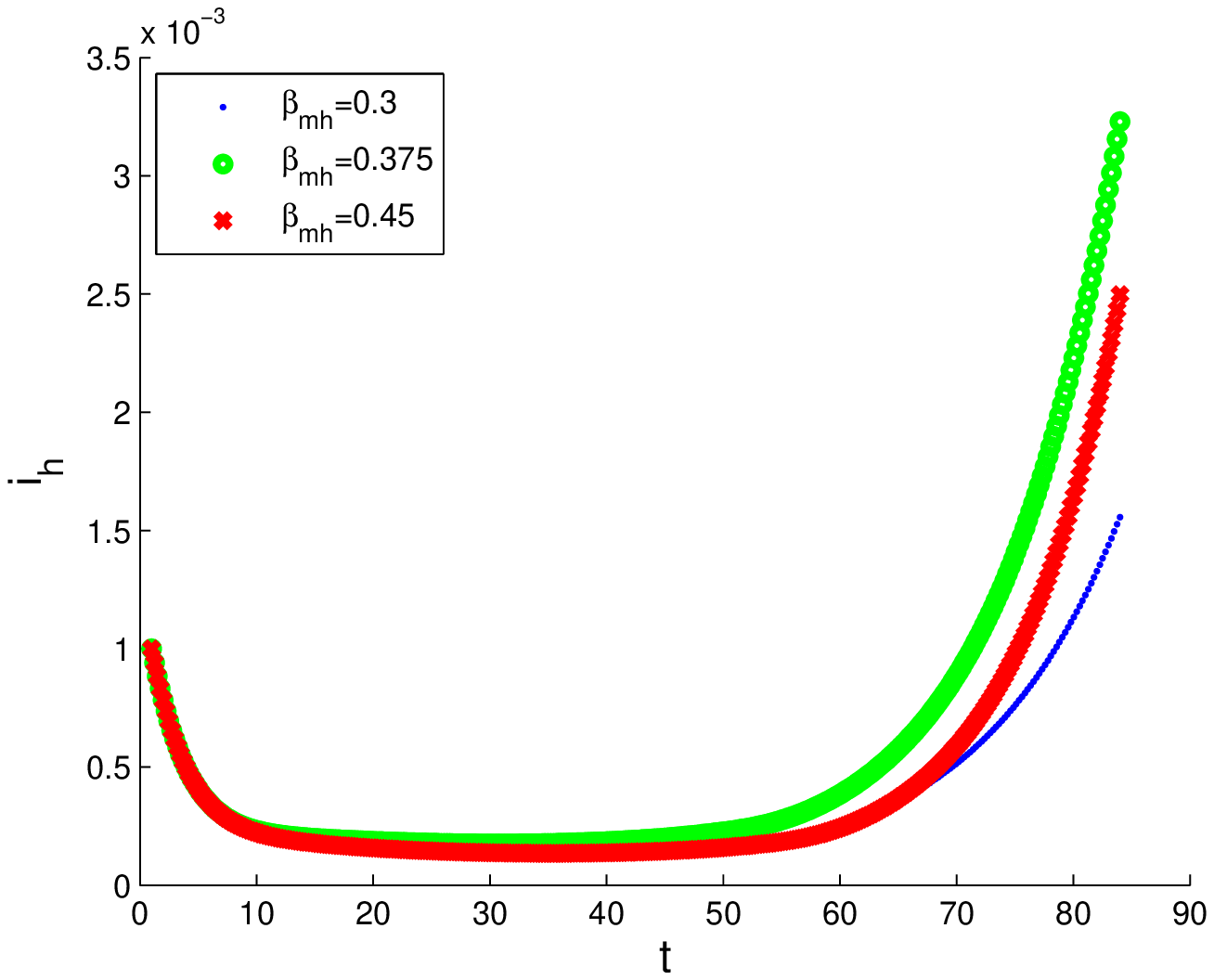,width=0.45\textwidth}\label{fig:dynamics:knee:4}}\\
\caption{Dynamics of control and state variables for varying values of
$\beta_{hm}$ and $\beta_{mh}$, corresponding to knee solution.}\label{fig:dynamics:knee}
\end{figure}

Figure~\ref{fig:dynamics:knee} shows the dynamics of the state and control variables
for different values of the transition probabilities, corresponding to the knee solution.
It can be seen that more controls are required for higher values of $\beta_{hm}$
and $\beta_{mh}$ during the entire period of study, which is consistent
with expectations. However, a higher value of  transition probability does not
necessarily mean a higher value of infected human population for the knee solution.
This is because the change in the shape of the Pareto front caused
by varying a value of transition probability.

\begin{figure}
\centering
\epsfig{file=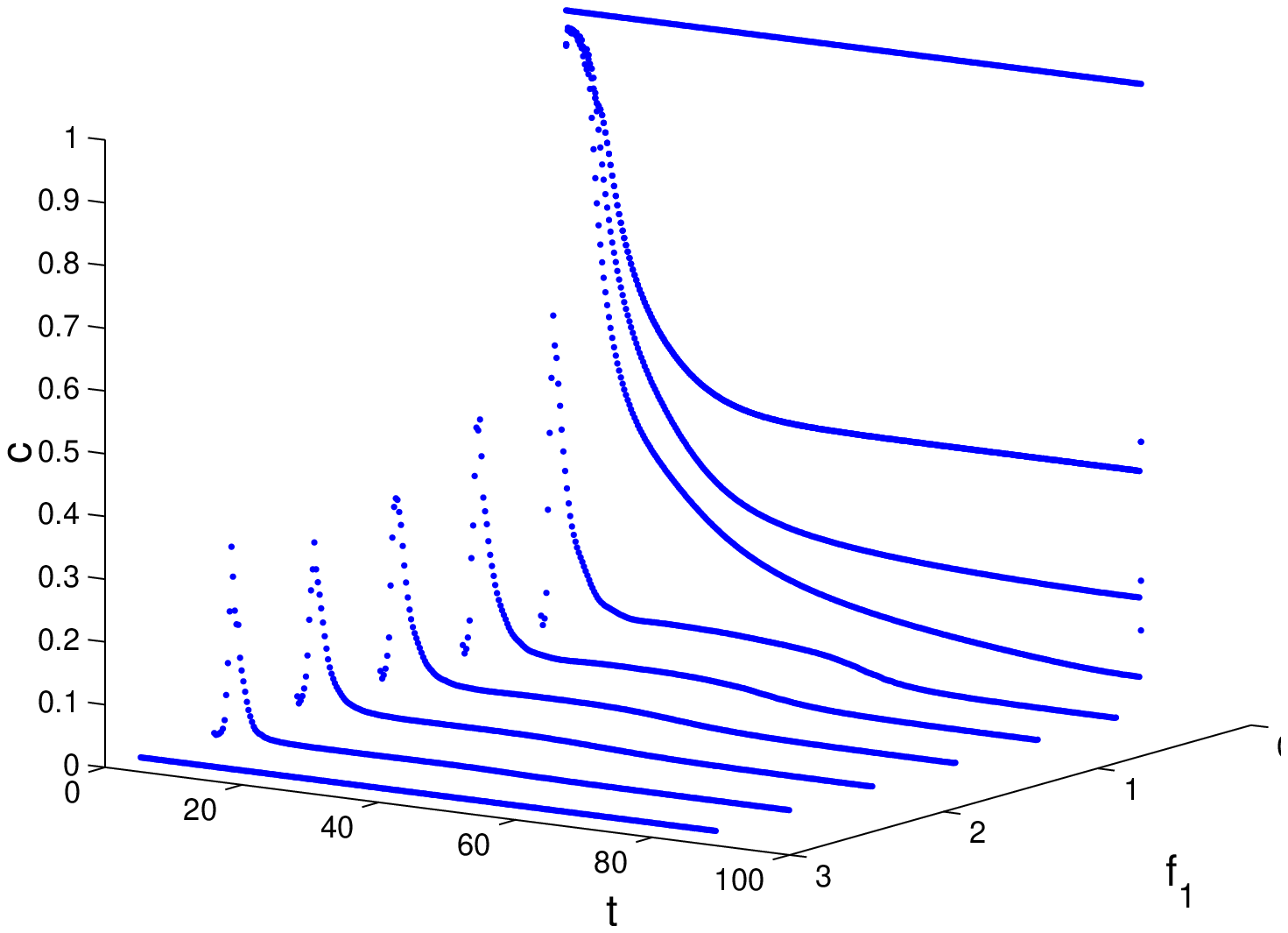,width=0.6\textwidth}
\caption{Discrete representation of the Pareto set.}\label{fig:suface:c}
\end{figure}
\begin{figure}
\centering
\subfloat[$0 \leq i_h(t) \leq 0.001$]{%
\epsfig{file=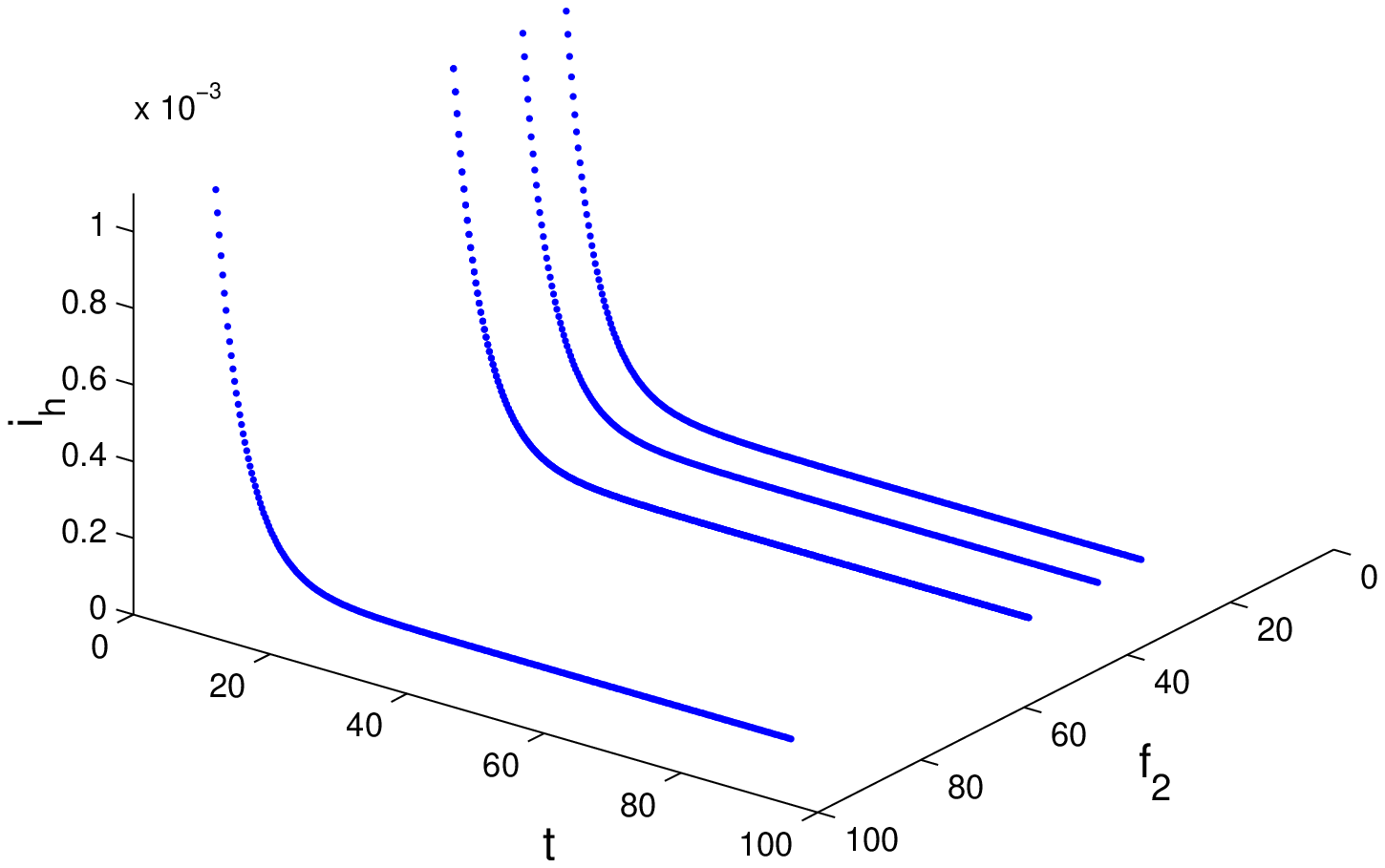,width=0.45\textwidth}\label{fig:suface:ih:1}}%
\subfloat[$0.01 \leq i_h(t) \leq 0.8$]{%
\epsfig{file=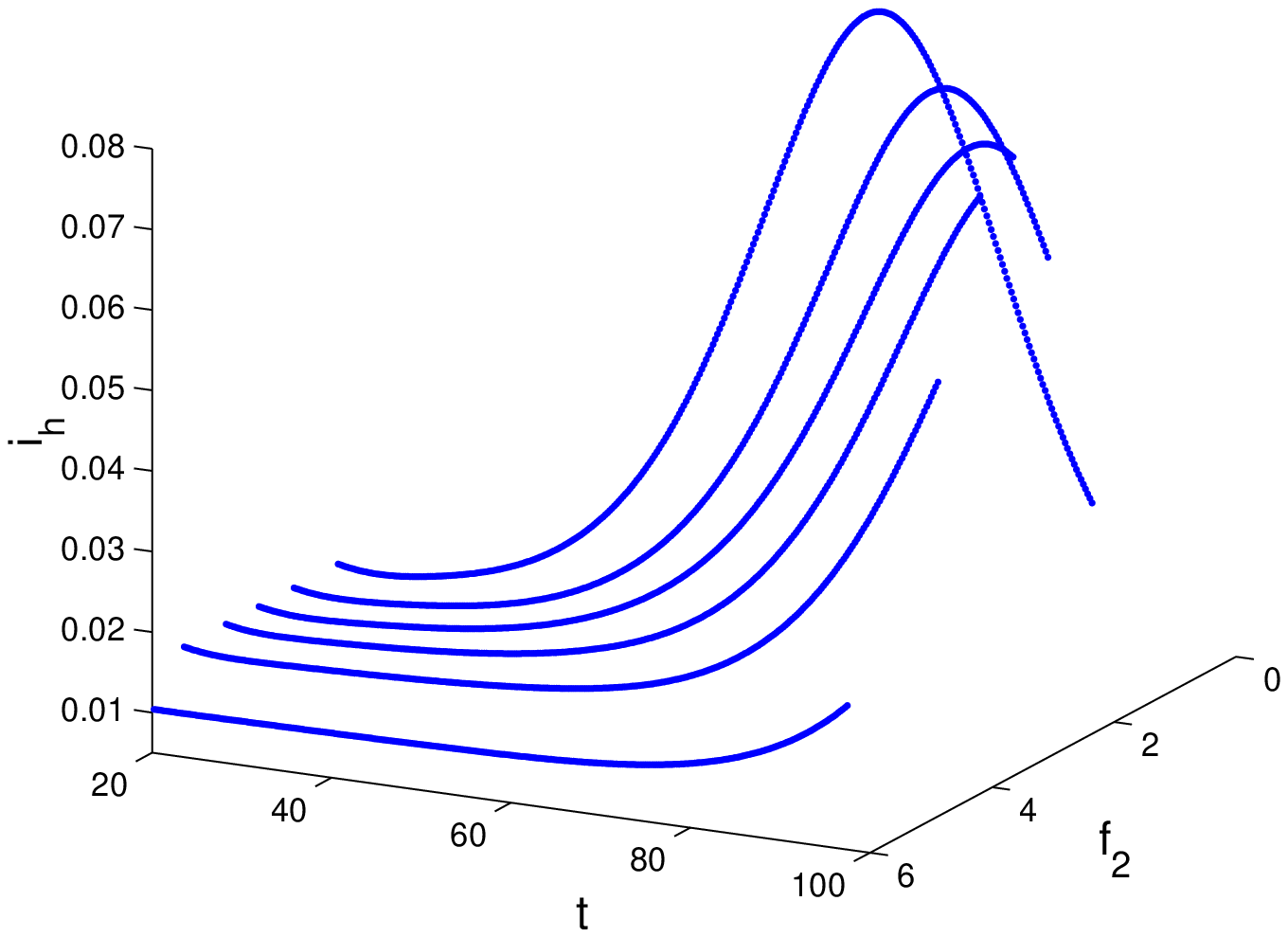,width=0.45\textwidth}\label{fig:suface:ih:2}}%
\caption{Discrete representation of infected humans.}\label{fig:suface:ih}
\end{figure}

From the Karush-Kuhn-Tucker condition, it can be induced that under smoothness
assumptions the Pareto set of a continuous multiobjective optimization problem
defines a piecewise continuous $(m-1)$-dimensional manifold. The Pareto set of
a continuous biobjective optimization problem is a piecewise continuous curve
in $\mathbb{R}^n$~\cite{ZhangZhou2008}. The optimal control for~\eqref{dengue:mop}
can be defined as a surface in $i_h(t)\times t$. Solving~\eqref{dengue:mop}
gives a discrete representation of this surface. Thus, Figure~\ref{fig:suface:c}
shows a discrete representation of a surface defining the Pareto set
for \eqref{dengue:mop}. For a given value of $f_1$, slicing the surface gives
the optimal control trajectory. From the figure, it can bee seen how this dynamic
changes from $c(\cdot) \equiv 0$ to $c(\cdot) \equiv 1$. A peak is observed in an early period
of $T$, which grows in accordance with the decrease in the number of infected humans.
Similarly, Figure~\ref{fig:suface:ih} shows a discrete representation of a surface
defining dynamics of infected humans across the Pareto set. A negative correlation
between the optimal control and infected population can be observed. As the amount
of control decreases, the number of infected humans increases. For higher values
of the control, the peak in $i_h(t)$ is smaller and occurs later. When the control
is reduced, the peak of disease outbreak grows up taking place earlier.

From the above discussion, it can be seen that the results for the optimal control problem,
presented in Section~\ref{sec:problem}, that are obtained using a  multiobjective
optimization approach provide comprehensive insights about optimal strategies
for dealing with the dengue epidemic and dynamics resulting from implementing
those strategies. The ability of trade-off solutions to reflect the underlying
nature of the problem constitutes a major advantage of multiobjective optimization,
motivating its practical use in the process of planning intervention measures
by health authorities.


\section{Conclusions}
\label{sec:cons}

Due to difficulties in treating the dengue disease, controlling and preventing
its outbreaks is essential for keeping people healthy, especially in regions
where the threat of dengue is high. This study discussed a mathematical model
for the dengue disease transmission from the optimal control point of view.
Multiobjective optimization approach is suggested for finding the optimal
control strategies to deal with an outbreak of the dengue epidemic.
A biobjective optimization problem is formulated, involving minimization
of expenses due to the infected population and costs of applying insecticide.
The approach avoids the use of \emph{a priori} information provided by the
decision maker in the form of weight coefficients and allows to reflect
the intrinsic nature of the problem.

The problem is numerically solved by discretizing the control variable
and using scalarization methods for approximating multiple Pareto optimal solutions.
The obtained results reveal different perspectives on applying insecticide:
a low number of infected humans can be achieved spending larger financial resources,
whereas low spendings for prevention campaigns result in significant portions
of the population affected by the disease. Different trade-offs between
the objectives are represented by the obtained solutions. Varying transmission
probabilities introduces changes in relation between the objectives and their trade-offs.
Because the problem is convex, the Pareto set defining the optimal control remains
unchanged for different values of transmission probabilities. The results of the
study suggest advantages of multiobjective optimization for finding the optimal control.
Once the Pareto set is approximated, the final decision on the control strategy
can be made taking into consideration available financial resources
and goals of public health care.

As future work, it is intended to extend the model for the dengue transmission
including different control variables that represent distinct measures for
controlling the disease. Investigating impacts of different values of model
parameters is also the subject of future work.


\section*{Acknowledgments}

This work has been supported by FCT (Funda\c{c}\~{a}o para a Ci\^{e}ncia e Tecnologia)
in the scope of projects UID/CEC/00319/2013 (ALGORITMI R\&D Center) and
UID/MAT/04106/2013 (Center for Research and Development in Mathematics
and Applications -- CIDMA). Torres was also supported by the project
PTDC/EEI-AUT/1450/2012, co-financed by FEDER under POFC-QREN
with COMPETE reference FCOMP-01-0124-FEDER-028894.
The authors are grateful to two referees
for valuable comments and suggestions.



\end{document}